 \documentclass{amsart}
  \usepackage {graphicx}
  \usepackage{geometry}
 \usepackage{amssymb}
 \usepackage{amscd}
 \usepackage [cmtip, all] {xy}
 \usepackage{hyperref}

 \newtheorem {theorem} {Theorem} [section]
\newtheorem {lemma} [theorem] {Lemma}

\newtheorem {corollary} [theorem]{Corollary}

\newtheorem {proposition} [theorem]{Proposition}

\theoremstyle {definition}

\newtheorem* {definition*}{Definition}
\newtheorem {example} [theorem] {Example}
\newtheorem {conjecture} [theorem] {Conjecture}
\newtheorem {problem} [theorem] {Problem}

\theoremstyle {remark}

\numberwithin {equation} {section}

 \def\sk{\operatornamewithlimits{sk}}
 \def\R{{\mathbb{R}}}
 
 \def\Z{{\mathbb{Z}}}

 \def\Cl{\operatorname{Cl}}

\begin{document}

\title{On approximability by embeddings
 of cycles in the plane}
 \author{ Mikhail Skopenkov}
 \address{Dept. of Diff. Geometry, Faculty of Mechanics and Mathematics,
 Moscow State University, Moscow, 119992, Russia,
 and Independent University of Moscow, B. Vlasyevsky, 11, 119002, Moscow, Russia.}
 \email{skopenkov@rambler.ru}
 \thanks{The author was supported in part by INTAS grant 06-1000014-6277,
Russian Foundation of Basic Research grants 05-01-00993-a,
06-01-72551-NCNIL-a, 07-01-00648-a, President of the Russian
Federation grant NSh-4578.2006.1, Agency for Education and Science
grant RNP-2.1.1.7988, and Moebius Contest Foundation for Young Scientists.}
 \subjclass{Primary:  57Q35; Secondary 54C25, 57M20}
 \keywords{approximability by embeddings, the van Kampen obstruction,
 line graph, derivative of a graph, derivative of a simplicial map, operation $d$,
 transversal self-intersection, standard $d$-winding, simplicial map,
 thickening}
 \begin{abstract} We obtain a criterion for
 approximability  of  piecewise linear  maps  $S^1\to\R^2$ by embeddings,
 analogous to the one proved by Minc for piecewise linear maps $I\to\R^2$.

 {\bf Theorem.} {\it Let $\varphi:S^1\to \R^2$ be a piecewise linear map, which is simplicial
 for some triangulation of $S^1$ with $k$ vertices. The map $\varphi$ is
 approximable by embeddings if and only if for each $i=0,\dots,k$ the
 $i$-th derivative $\varphi^{(i)}$ (defined by Minc) neither
 contains transversal self-intersections nor is the standard winding of
 degree $\not\in\{-1,0,1\}$.}

 We deduce from the Minc result
 the completeness of the van Kampen obstruction to approximability by
 embeddings of piecewise linear maps $I\to\R^2$.
 We also generalize these criteria to simplicial maps $T\to S^1\subset\R^2$,
 where $T$ is a graph without vertices of degree $>3$.
 \end{abstract}
 \maketitle

 \section{Introduction}

 A PL map $\varphi:K\to \R^2$ of a graph $K$ is {\it approximable by
 embeddings} in the plane, if for each $\varepsilon>0$ there is an
 $\varepsilon$-close to $\varphi$ map $f:K\to\R^2$ without
 self-intersections. In the major part of this paper we consider
 the case when $\varphi$ is either a path or a cycle, i. e.
 either $K\cong I$ or $K\cong S^1$.

 \begin{example} \cite{11}
 The standard $d$-winding $S^1\to S^1\subset\R^2$ is
 approximable by embeddings in the plane if and only if
 $d\in\{-1,0,1\}$.
 \end{example}

 It can be also proved that a simplicial map $S^1\to S^1$ is
 approximable by embeddings if and only if its degree $d\in\{-1,0,1\}$
 (see Theorem~1.3).
 A {\it transversal self-intersection} of a PL map $\varphi:K\to \R^2$
 is a pair of disjoint arcs $i,j\subset K$ such that
 $\varphi i$ and $\varphi j$ intersect transversally in the plane.

 \begin{example} An Euler path or cycle in
 a graph in the plane
 is approximable by embeddings if and only if it
 does not have transversal self-intersections
 (hence any Euler graph in the plane has an Euler cycle,
 approximable by embeddings).
 \end{example}

 The notion of approximability by embeddings appeared
 in studies of embeddability of compacta into $\R^2$
 (see \cite{11, 14, 12}, for recent surveys see
 \cite{17}, \cite[\S9]{7}, \cite[\S4]{2}, \cite[\S1]{8},
 we return to this topic in the end of \S1.)
 There exists an algorithm of checking
 whether a given simplicial map is
 approximable by embeddings (see \cite{13},
 or else Simple-minded Criterion~4.1  below).
 A more convenient to apply criterion for  approximability by embeddings
 of a simplicial path in the plane was proved in
 \cite{6} (Theorem~1.3.I below, generalizing Example~1.2).  The main
 result of this paper is an analogous criterion for approximability
 by embeddings of a cycle in the plane (Theorem~1.3.S below, also
 generalizing Example~1.2). These criteria assert that, in some sense,
 transversal self-intersections are the only obstructions
 to approximability by embeddings.
 Clearly, this is not true literally \cite{11},
 and there is no Kuratowsky-type criterion.

 We state our criterion (Theorem~1.3) in terms of {\it the derivative}
 of a path \cite{5}, \cite[''{\it the operation $d$}'']{6}.
 Let us give the definition
 (Fig.~1).
 First let us define {\it the derivative} $G'$ of a graph $G$
 (also called {\it adjoint, conjugate, covering, derived, edge, edge-to-vertex dual, interchange, line, representative graph}).
 The vertex set of the graph $G'$ is in 1-1 correspondence
 with the edge set of $G$.
 For an edge $a\subset G$ denote by $a'\in G'$ the corresponding vertex.
 Vertices $a'$ and $b'$ of $G'$
 are joined by an edge if and only if the edges $a$ and $b$ are adjacent
 in $G$.
 Note that the derivatives
 $G'$ and $H'$ of homeomorphic but not isomorphic graphs $G$
 and $H$ are not necessarily homeomorphic. Beware that we slightly modify the definition of $G'$ below.

 Now let $\varphi$ be a path in the graph $G$ given by a sequence
 of vertices $x_1,\dots,x_k\in G$, where
 $x_i$ and $x_{i+1}$ are joined by an edge.
 Then $(x_1x_2)',\dots,(x_{k-1}x_k)'$
 is a sequence of vertices of $G'$.
 In this sequence replace each segment
 $(x_ix_{i+1})'$, $(x_{i+1}x_{i+2})$, $\dots$, $(x_{j-1}x_{j})'$ such that
 $(x_ix_{i+1})'=(x_{i+1}x_{i+2})'=\dots=(x_{j-1}x_{j})'$ by a single
 vertex.
 The obtained sequence of vertices determines a path in the graph $G'$.
 This path $\varphi'$ is called {\it the pre-derivative} of the path $\varphi$. We consider both paths $\varphi$ and $\varphi'$ as piecewise-linear maps  $\varphi\colon I\to G$ and $\varphi'\colon I\to G'$.

 \begin {figure}
\includegraphics {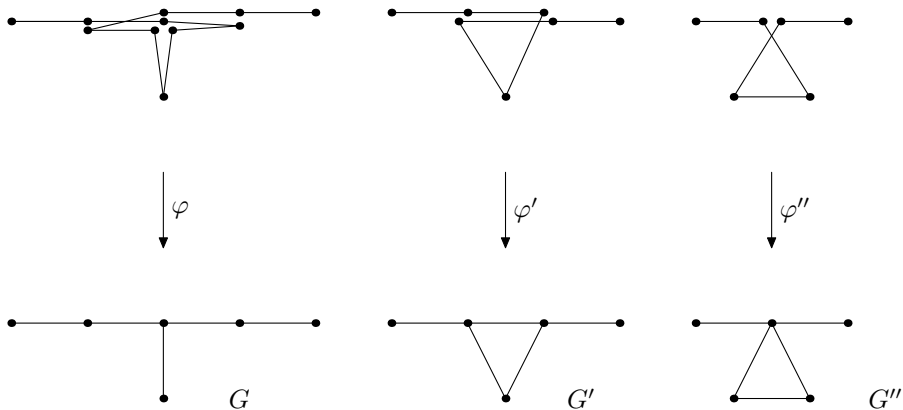}
\caption {Derivatives of graphs and paths}
\label {t-fig1}
\end {figure}

 Now let $G$ be embedded into the plane so that the edges are rectilinear segments. The derivative of $G$ might not be planar; e.g., the derivative of a $5$-od (the cone over $5$ points) is the nonplanar Kuratowsky graph.
 But if the path
 $\varphi$ does not have transversal self-intersections,
 then the image of the map $\varphi'\colon I\to G'$ is a planar subgraph
 $G'_\varphi\subset G'$; we give the construction of a natural
 embedding $G'_\varphi\to\R^2$ in \S2, Definition~of~$N'$.
 Change $G'$ to the image $G'_\varphi$ and $\varphi'$ to its onto
 restriction $\varphi':I\to G'_\varphi$ called {\it the derivative} of $\varphi$. In what follows $ G' $ and $\varphi '$ denote these restrictions.
 Define the $k$-th \emph{derivative}
 $\varphi^{(k)}$ inductively. Surely, the map $\varphi ^ {(k)}\colon I\to\mathbb{R}^2$ is well-defined, only if all the previous derivatives $\varphi, \varphi',\dots, \varphi ^ {(k-1)}\colon I\to\mathbb{R}^2$ have no  transversal self-intersections. For a cycle $\varphi$ the definition of
 {\it the derivative} cycle $\varphi'$ is analogous (the domain of $\varphi'$ might degenerate to a point or an empty set).

 An example to be
 used in the sequel is that $\varphi'=\varphi$ for a standard
 $d$-winding $\varphi:S^1 \to S^1$ with $d\ne0$.  Clearly,
 $\varphi'$ is an embedding for any Euler path or cycle $\varphi$ without transversal self-intersections.
 Thus Example~1.2 is indeed a specific case of the following theorem.

 \begin{theorem} I) \cite{6}
 Let $\varphi:I\to \R^2$ be a PL map, which is simplicial for some
 triangulation of $I$ with $k$ vertices. The map $\varphi$ is
 approximable by embeddings if and only if for each $i=0,\dots,k$
 the $i$-th derivative $\varphi^{(i)}$ does not contain
 transversal self-intersections.

 S) Let $\varphi:S^1\to \R^2$ be a PL map, which is simplicial for some
 triangulation of $S^1$ with $k$ vertices. The map $\varphi$ is
 approximable by embeddings if and only if for each
 $i=0,\dots,k$ the $i$-th derivative $\varphi^{(i)}$
 neither contains transversal self-intersections
 nor is the standard winding of degree $d\not\in\{-1,0,1\}$.
 \end{theorem}

 Let us comment the theorem. The maps $\varphi, \varphi^{(1)},\varphi^{(2)},\dots \colon I\to\mathbb{R}^2$ or $S^1\to\mathbb{R}^2$ are well-defined as maps into the plane up to the first number for which the derivative has a transversal self-intersection. If such a number exists, then the condition after `if and only if' already fails, and $\varphi$ is not approximable by embeddings. Otherwise, for 1.3.I, the condition holds, and $\varphi$ is approximable by embeddings (whereas for 1.3.S one also needs to find the degree of the winding $\varphi ^ {(k)}$; cf.~Lemma~2.3).

 We prove both 1.3.I and 1.3.S in \S2.
 Our proof of 1.3.I is simpler than the one given in \cite{6}.

 In \S3 we apply Theorem~1.3 to prove the following criterion.

 \begin{corollary}
 A PL map $\varphi:I\to \R^2$ is approximable by
 embeddings if and only if one of the following equivalent
 conditions holds:

 D) (the deleted product property)
 There is a map $\{ (x,y)\in I\times I\: x\ne y\}\to S^1$ such
 that its restriction to the set
 $\{ (x,y)\in I\times I: \varphi x\ne\varphi y\}$ is
 homotopic to the map given by the formula
 $\tilde\varphi(x,y)=\frac {\varphi x-\varphi y} {\|\varphi x-\varphi y\|}$;

 V) the van Kampen obstruction (defined in \S3) $v(\varphi)=0$.
 \end{corollary}

 The criterion 1.4.V, although more difficult to state,
 is more easy to apply than 1.3.I and 1.4.D.
 In Corollary~1.4 the arc $I$ cannot be replaced by $S^1$: the standard
 3-winding is a counterexample \cite{8}.
 Obstructions like 1.4.D and 1.4.V appear in the related theory
 of approximability by link maps, i. e. by maps with
 disjoint images, but the criteria analogous to 1.3.I and 1.4.DV
 are not true (Example~3.3 below).
 In \S3 we show that the $n$-dimensional generalizations of
 conditions 1.4.D and 1.4.V are equivalent for any PL maps
 $\varphi:K^n\to \R^{2n}$ (Proposition~3.2 below).

 In \S4 we generalize criteria~1.3 and~1.4 to PL maps
 $\varphi:K\to \R^2$, where $K$ is an arbitrary graph.
 We prove the following theorem (see Definition of the derivative in \S2).

 \begin{theorem}
 Let $T$ be a graph without vertices of degree $>3$.
 Suppose that $T$ has $k$ vertices.
 A simplicial map $\varphi:T\to S^1\subset\R^2$ is approximable
 by embeddings if and only if the van Kampen obstruction
 $v(\varphi)=0$ and $\varphi^{(k)}$ does
 not contain standard windings of degree $d\ne\pm 1$, $d$ odd.
 \end{theorem}

 \begin{conjecture} Theorem~1.5 is true
 for a simplicial map $\varphi:K\to G\subset\R^2$, where $K$ and $G$
 are arbitrary graphs.
 \end{conjecture}

 If Conjecture~1.6 is true, then
 {\it a simplicial map $\varphi:T\to\R^2$ of a tree $T$
 is approximable by embeddings if and only if $v(\varphi)=0$}
 \cite[Problem~4.5]{2}.

 \begin{conjecture} A piecewise linear path $\varphi:I\to \mathbb{R}^2$ is approximable by embeddings
 if and only if for each pair of arcs $I_1,I_2\subset I$ such that $I_1\cap I_2=\emptyset$ the pair of restrictions
 $\varphi:I_1\to \mathbb{R}^2$ and $\varphi:I_2\to \mathbb{R}^2$ is approximable by link maps (i. e. maps with disjoint images).
 \end{conjecture}

 We conclude \S1 by some words on the history of the
 notion of approximability by embeddings.
 We define the decomposition of a 1-dimensional compactum
 into {\it an inverse limit} and show how the notion
 of approximability by embeddings appears in studies of planarity
 of this compactum.
 We do not use this definition in our paper.
 To give an example,
 let us construct {\it the 2-adic van Danzig solenoid}.
 Take a solid torus $T_1\subset\R^3$.
 Let $T_2\subset T_1$ be a solid torus going twice along the axis of
 the torus $T_1$.
 Analogously, take $T_3\subset T_2$ going twice
 along the axis of $T_2$.
 Continuing in the similar way, we obtain an infinite sequence
 of solid tori $T_1\supset T_2\supset T_3\supset\dots$
 The intersection of all tori $T_i$ is a 1-dimensional compactum and is
 called  {\it the 2-adic van Danzig solenoid}.
 By {\it the inverse limit} of an infinite sequence
 of graphs and simplicial maps between them
 $
 \xymatrix@1{
 {K_1} & {K_2}\ar[l]_{\varphi_1} & {K_3} \ar[l]_{\varphi_2} & {\dots} \ar[l]_{\varphi_3}
 }
 $
  we mean the compactum
 $$
 C = \{\, (x_1, x_2, \dots) \in l_2 \,
 : \, x_i\in K_i \text { and } \varphi_ix _ {i+1} =x _ {i} \, \}.
 $$
 One can see from our construction that
 for the van Danzig solenoid all $K_i\cong S^1$
 and all $\varphi_i$ are 2-windings.
 It can be proved that any 1-dimensional compactum can be represented
 as an inverse limit.
 Such representation shows that any 1-dimensional compactum can be embedded
 into $\R^3$. It also gives an easy sufficient condition to planarity:
 for each positive integer $i$
 there should exist an embedding $f_i:K_i\to \R^2$  such that the map
 $f_i\circ \varphi_i$ is approximable by embeddings and
 $f_{i+1}$ is $2^{-i}$-close to $f_i\circ\varphi_i$.

 \section{Proofs}

 Theorem~1.3 follows from Example 1.1
 and Lemmas 2.1, 2.2.A (for $K\cong I, S^1$)
 and 2.3, which are interesting results in themselves.

 \begin{lemma} (for $K\cong I$ see \cite{6})
 Suppose that a simplicial map $\varphi:K\to G\subset\R^2$
 of a graph $K\cong S^1$ or $K\cong I$ does not have
 transversal self-intersections.
 If $\varphi'$ is approximable by
 embeddings, then $\varphi$ is approximable by embeddings.
 \end{lemma}

 \begin{lemma} A) \cite{6}
 If a simplicial map $\varphi:K\to G\subset\R^2$ is approximable by
 embeddings, then the map $\varphi'$ is approximable
 by embeddings.

 V) If a simplicial map $\varphi:K\to G\subset\R^2$ is approximable by
 $\mod 2$-embeddings, then the map $\varphi'$ is approximable by
 $\mod 2$-embeddings.
 \end{lemma}

 Here {\it a $\mod2$-embedding} is a general position map $f:K\to\R^2$
 such that for each pair $a,b$ of disjoint edges of $K$ the set
 $fa\cap fb$ consists of an even number of points.
 Definition of the derivative $\varphi'$ needed for Lemma~2.2 is presented below

 \begin{lemma} Let $\varphi:S^1\to G$ be a PL map,
 which is simplicial for some
 triangulation of $S^1$ with $k$ vertices.
 Then either the domain of $\varphi^{(k)}$ is empty or $\varphi^{(k)}$ is
 a standard winding of degree $d\ne0$.
 \end{lemma}

 This number $d$ can be considered as the generalization
 of {\it the degree} of any simplicial map $S^1\to G$.
 So it is interesting to get the solution of the following problem
 (it may also make criteria~1.3 and~1.5 more easy to apply).

 \begin{problem}
 Find an easy algorithm for calculation of the degree
 of the winding $\varphi^{(\infty)}$ for a given PL map $\varphi:S^1\to G$.
 \end{problem}

 Futher we use the following generalization of the definition
 of the derivative of a path stated in \S1.

 \begin{definition*}[Definition of the derivative \cite{6}, see Fig.~1 and a part of
 Fig.~4]  First let us construct the graph $K'_\varphi$, which is the
 domain of the derivative $\varphi'$.  By {\it a $\varphi$-component}
 of the graph $K$ we mean any connected component $\alpha$ of
 $\varphi^{-1}a$ mapped {\it onto} $a$, for some edge $a\subset G$.
 The vertex set of $K'_\varphi$ is in 1-1 correspondence
 with the set of all $\varphi$-components.
 For a $\varphi$-component $\alpha\subset K$ denote by
 $\alpha'\in K'_\varphi$ the corresponding vertex.
 Vertices $\alpha'$ and $\beta'$ are joined by an edge in $K'_\varphi$
 if and only if $\alpha\cap \beta\ne\emptyset$.
 The {\it derivative} $\varphi': K'_\varphi\to G'$
 is a simplicial map defined on the vertices $K'_\varphi$
 by the formula $\varphi'\alpha'=(\varphi \alpha)'$.
 Change $\varphi'$ to its onto restriction
 $\varphi':K'_\varphi\to \varphi'K'_\varphi$.
 (In the original definition
 \cite{6} $G'$ is denoted by $D(G)$,
 $\varphi'$ by $d[\varphi]$ and $K_\varphi'$ by $D(\varphi,K)$.)
 \end{definition*}

 \begin{proof}[Proof of 2.3]
 We say that a simplicial map $\varphi:K\to G$ is
 {\it ultra-nondegenerate}, if for each edge $a\subset K$
 the image $\varphi a$ is an edge of $G$ and
 for each pair $a,b\subset K$
 of adjacent edges we have $\varphi a\ne\varphi b$.
 Denote by $|K|$ the number of vertices in a graph $K$.
 Clearly, if $K\cong S^1$,
 then $|K_\varphi'|\le |K|$, and $|K_\varphi'|=|K|$ only if
 $\varphi$ is ultra-nondegenerate.
 Therefore it suffices to prove the lemma for this latter case
 (because the cases $K_\varphi'\cong I$ or $K_\varphi'$ is a point
 are trivial).
 In this case the lemma is obvious, but we give the proof.

 Let us prove that if an ultra-light simplicial onto map $\varphi:K\to G$
 of the graph $K\cong S^1$ is not a standard winding of a nonzero degree,
 then $|G'|>|G|$.
 Note that for ultra-nondegenerate $\varphi:S^1\to G$ the graph $G$
 does not contain hanging vertices.
 If the degree of each vertex of $G$ is two,
 then $\varphi$ is an ultra-nondegenerate simplicial map $S^1\to S^1$,
 consequently $\varphi$ is a standard winding, that contradicts
 to our assumption.
 So $G$ contains a vertex of degree at least 3.
 Then by the above the number of edges of $G$ is greater than
 the number of vertices, hence $|G'|>|G|$.
 Since for a simplicial onto map $\varphi:K\to G$ we have
 $1\le |G|\le |K|$, it follows that $|G|,|G'|,\dots,|G^{(k)}|\le k$
 (recall that we define $\varphi'$ to be an onto map).
 This yields that one (and then $k$-th) of the derivatives
 $\varphi,\dots,\varphi^{(k)}$ is a standard winding of a nonzero degree,
 because otherwise we obtain $1+k \le |G|+k \le |G^{(k)}| \le k$.
 \end{proof}

 Now let us give the proposed construction of the
 embedding $G'_\varphi\to\R^2$.  It is more convenient for us to
 consider {\it thickenings} of the graphs rather than embeddings of the
 graphs into the plane.  Then the proposed construction is equivalent
 to the construction of {\it the derivative of a thickening}
 (Definition~of~$N'$ below).  Further we assume that {\it a thickening
 $N$} of the graph $G$ in the plane (i. e., a regular neighbourhood of
 $G\subset\R^2$) is fixed.  We also assume that a handle decomposition
 (denoted by $S$)
 $$
 N=\bigcup\limits_{x\in\text{vertex set of }G}N_x
    \cup\bigcup\limits_{a\in\text{edge set of }G}N_{(a)}
 $$
 corresponding to the graph $G$ is also fixed,
 where $N_x$ are 2-discs and
 $N_{(a)}$ are joining them strips.
 Denote by
 $N_a$ the restriction $N_x\cup N_{(a)}\cup N_y$
 of $N$ to an edge $a=xy$.
 Actually, we do not use the planarity of $N$ in our proofs, the thickening
 $N$ can be assumed to be just orientable (orientability is needed for
 Example~1.1).
 Let us state the definition of the derivative $N'$ of a thickening $N$.
 This thickening $N'$
 depends on the simplicial map $\varphi:K\to G\subset N$
 and is well-defined only if $\varphi$ does not contain transversal
 self-intersections. Moreover,
 for an arbitrary $K$ we must also assume
 that there are no pairs of arcs $i,j\subset K$
 (not necessarily disjoint!) such that the intersection
 $\varphi i\cap\varphi j$ is transversal.

 \begin{definition*}[Definition~of~$N'$, see Fig.~2]
 Let $\varphi:K\to G\subset N$ be a simplicial
 map such that for any pair of arcs $i,j\subset K$
 the intersection $\varphi i\cap\varphi j$ (maybe empty) is
 not transversal.
 Let us construct discs $N'_{a'}$
 for each vertex $a'\in G'$ and strips $N'_{(a'b')}$
 for each edge $a'b'\subset G'$.
 Then $N'$ together with its handle decomposition $S'$
 is defined by the formula
 $N'=\bigcup N'_{a'}\cup\bigcup N'_{(a'b')}$.
 Here we take $N'_{a'}=N_{(a)}$ for each edge $a\subset G$.
 For each pair $a,b\subset G$ of adjacent edges such that
 $(\varphi')^{-1}(a'b')\ne\emptyset$
 we join the two discs $N'_{a'}$ and $N'_{b'}$ by a narrow strip
 $N'_{(a'b')}$ in $N_{a\cap b}$.
 Since the intersection of arcs $a\cup b$ and
 $c \cup d$ is not transversal for any pair of adjacent edges
 $c,d\subset K$,
 it follows that we can choose the strips $N'_{(a'b')}$
 so that they do not intersect for distinct $a'b'$.
 \end{definition*}

 \begin {figure}
\includegraphics{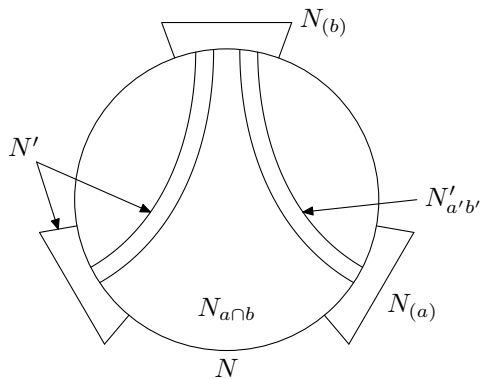}
\caption {Derivative of a thickening}
\label {t-fig2}
\end {figure}

 This definition can also be considered
 as a construction of an embedding $N'\to N$, and also
 $G_\varphi'\to\R^2$.
 Note that $S'$ and the topological type of $N'$ do not depend on the
 choice of the strips $N'_{(a'b')}$ in our definition.
 The alternative definition of the derivative thickening $D(N)$ in \cite{6}
 does not depend also on the map $\varphi$.
 The thickening
 $N'$ of our paper means the subthickening of $D(N)$ of \cite{6},
 corresponding to the subgraph $G_\varphi'\subset G'$.

 Clearly, for investigation of approximability by
 embeddings of simplicial maps $K\to G\subset\R^2$
 it suffices to consider only the approximations $f:K\to N$.
 Now we are going to reduce the problem of approximability
 by embeddings of a given map to the problem of existance of
 an embedding close to it in some sense ($S$-close to it).

 \begin{definition*}[Definition~of~an~$S$-approximation, cf. \cite{6}]
 A map $f: K \to N$
 is {\it an $S$-approximation} of the map $\varphi$,
 or $f$ is {\it $S$-close} to $\varphi$,
 if the following conditions hold:
 \begin{enumerate}
 \item $fx\subset N_{\varphi x}$ for each vertex or edge $x$ of $K$
 \item $x\cap f^{-1}N_{(\varphi x)}$ is connected for each edge
       $x$ of $K$ with nondegenerate $\varphi x$.
 \end{enumerate}
 \end{definition*}

 Proposition~2.9 in \cite{6} asserts that the map
 $\varphi:K\to G$
 is approximable by embeddings if and only if there is an embedding
 $f:K\to N$, $S$-close to $\varphi$.

 A PL map $\varphi:K\to N$ is {\it degenerate},
 if $\varphi c$ is a point for some
 edge $c\subset K$.
 Now let us prove the following easy Contracting Edge Proposition~2.5
 that in some sense allows us to assume that in 2.1 and 2.2
 the map $\varphi$ is nondegenerate.

 \begin{proposition}[Contracting~Edge~Proposition] Let $\varphi:K\to G$
 be a simplicial map such that $\varphi c$ is a point for some
 edge $c\subset K$.
 Let $K/c$ be the graph obtained from $K$ by contracting
 the edge $c$, and let $\varphi/c:K/c\to G$ be the corresponding map.
 Then

 D) $K_{\varphi/c}'=K_\varphi'$, $G_\varphi'=G_{\varphi/c}'$
 and $(\varphi/c)'=\varphi'$.

 A) for $K\cong S^1$ or $K\cong I$ the map $\varphi/c$ is approximable
 by embeddings if and only if $\varphi$ is approximable by embeddings.

 K) for an arbitrary $K$ if $\varphi$ is approximable by embeddings,
 then $\varphi/c$ is approximable by embeddings.

 V) If $\varphi$ is approximable by $\mod2$-embeddings,
 then $\varphi/c$ is approximable by $\mod2$-embeddings.
 \end{proposition}

 \begin{proof}[Proof of 2.5]
 D) is obvious.

 A) Let us prove the direct implication.
 Let $f:K/c\to N$ be an embedding, $S$-close to $\varphi/c$.
 Let $a\subset K$ be an edge adjacent to $c$
 (if $c$ is a connected component of $K$, then the proposition is obvious).
 Add a new vertex to the edge $a$ of the graph $K/c$ (Fig.~3.a).
 Since $K\cong S^1$ or $K\cong I$, it follows that
 the obtained graph is isomorphic to $K$
 and the embedding $f:K\to N$ is the required.
 The reverse implication is a specific case of statement K).

 K) Let $f:K\to N$ be an embedding, $S$-close to $\varphi$.
 Make the move shown in Fig.~3.b.
 We obtain an embedding $\bar f:K/c\to N$, $S$-close to $\varphi/c$.

 V) Let $f$ be a $\mod2$-embedding, $S$-close to $\varphi$.
 Make the move shown in Fig~3.b.
 We obtain an $S$-close to $\varphi/c$ map $\bar f:K/c\to N$.
 It suffices to prove that $|\bar fa\cap \bar fb|=0\pmod2$
 for each pair of disjoint edges $a,b\subset (K/c)$.
 Indeed, both $a$ and $b$ are also edges of $K$,
 and at least one of them is not
 adjacent to $c$ (because $a$ and $b$ are disjoint in $K/c$).
 If neither $a$ nor $b$ is adjacent to $c$, then
 $|\bar fa\cap \bar fb|=|fa\cap fb|=0\pmod2$.
 If, for example, $b\in K$ is adjacent to $c$
 and $a$ is not adjacent to $c$, then
  $|\bar f a\cap \bar f b|=|fa\cap fb|+|fa\cap fc|=0\pmod2$,
 that proves the proposition.
 \end{proof}

 \begin {figure}
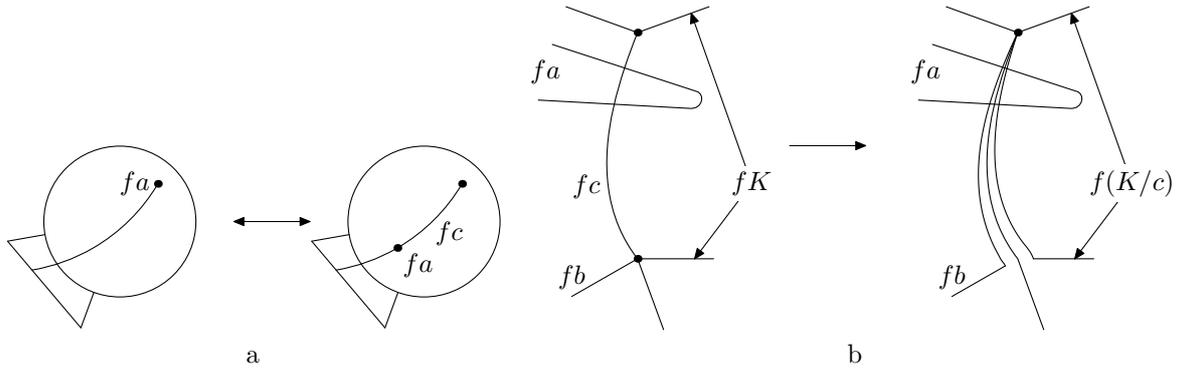

 \begin{tabular}{cc}
\includegraphics {skopenkov2.301} & \includegraphics{skopenkov2.302}\\
a & b
\end{tabular}
\caption {Moves for degenerate maps}
\label {t-fig3}
\end {figure}

 Degenerate maps appear in our proof of 2.1 and~2.2
 even if the map $\varphi:K\to G$ is nondegenerate.
 We are going to construct
 a graph $\bar K'_\varphi$ and a pair of (degenerate)
 simplicial maps $\xymatrix{ {G} & \bar K '_\varphi \ar[l]_{\bar\varphi} \ar[r]^{\bar\varphi'} & G' }$
 that can be obtained from $\varphi$ and $\varphi'$ respectively
 by the operation from Contracting~Edge~Proposition~2.5
 (this is true under some assumptions on $\varphi$,
 we present the details below).
 Together with the construction of the embedding $N'\to N$
 (see Definition of~$N'$ above) this immediately proves~2.1
 (Fig.~4, 5, 6).

 \begin{definition*}[Definition~of~$\bar\varphi$ and~$\bar\varphi'$, see Fig.~4]
 Suppose that the map $\varphi$ is nondegenerate and $K$ does not
 have vertices of degree 0.
 Take the disjoint union of all $\varphi$-components
 of $K$ (see Definition~of~$\varphi'$). Join by an edge any two vertices
 belonging to distinct $\varphi$-components and
 corresponding to the same vertex of $K$.
 Denote the obtained {\it semi-derivative} graph by $\bar K'_\varphi$.
 Thus a $\varphi$-component $\alpha\subset K$ is also a subgraph
 of $\bar K_\varphi'$ denoted by $\bar\alpha'$.
 Further we identify the points of $\alpha$ and $\bar\alpha'$.
 Let the simplicial maps $\bar\varphi$ and
 $\bar \varphi'$ be the evident projections
 $\bar K'_{\varphi}\to G$ and $\bar K_\varphi'\to G'$ respectively,
 defined on the vertex sets by $\bar\varphi x=\varphi x$ and
 $\bar \varphi'x=(\varphi \alpha)'$,
 where the vertex $x\in\bar K_\varphi'$ belongs to
 the $\varphi$-component $\bar\alpha'$.
 \end{definition*}

 \begin {figure}
\includegraphics {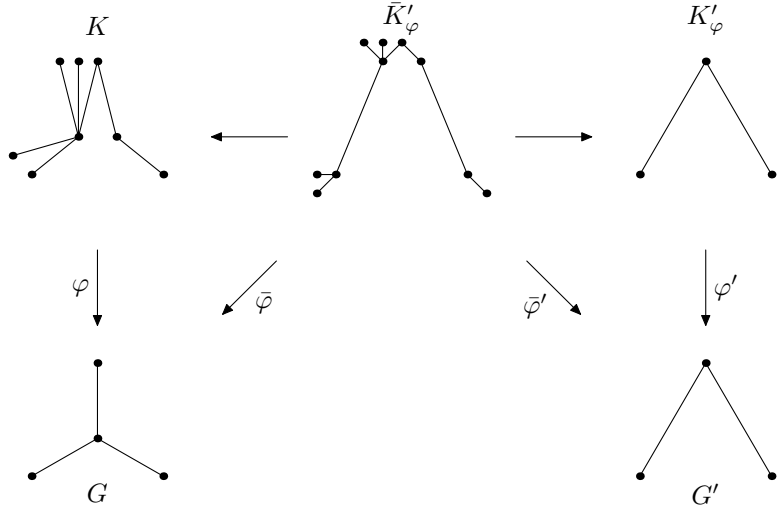}
\caption {Semi-derivatives of a simplicial map}
\label {t-fig4}
\end {figure}

 \begin{proof}[Proof of 2.1]
 By Contracting~Edge~Proposition~2.5.D,A the map $\varphi$ can be
 assumed to be nondegenerate. We also may assume that
 $K$ does not have vertices of degree 0.
 It can be easily checked that
 $\varphi$ and $\varphi'$ can be obtained from
 $\bar\varphi$ and certain restriction of $\bar\varphi'$
 respectively by the operation from
 Contracting Edge Proposition~2.5.
 If any two $\varphi$-components have at most
 one common point, then $\varphi'$
 can be obtained from $\bar\varphi$ itself in this way.
 But for $K\cong S^1$ this asumption is not satisfied only if
 $K$ has two $\varphi$-components. Evidently, the map $\varphi$
 is approximable by embeddings in this case.
 So it suffices to prove that

\smallskip
 (*) if $\bar\varphi'$ is approximable by embeddings,
 then $\bar\varphi$ is approximable by embeddings.
\smallskip

 We prove ($*$) for an arbitrary graph $K$. If $\bar \varphi'$ is approximable by embeddings,
 then there is a $S'$-close to $\bar \varphi'$
 embedding $\bar K_\varphi'\to N'$.
 Define the embedding $f:\bar K_\varphi'\to N$ to be the composition
 of this embedding and the embedding $N'\to N$
 constructed in Definition~of~$N'$ (Fig.~5, where this construction is applied
 to the map $\varphi$ from Fig.~4).
 Clearly, there exists a new handle decomposition $N=\bigcup \bar N_a \cup \bigcup \bar N_{(ab)}$,
 denoted by $\bar S$, such that $f$ is an $\bar S$-approximation of
 $\bar\varphi$ (Fig.~6, cf. \cite[Proposition~4.9]{6},
 or see some generalization of the decomposition, constructed
 in the proof of Lemma~4.5.D.)
 Then $f:\bar K'_{\varphi}\to \bar N$ (where $\bar N$ is
 $N$ with the new handle decomposition $\bar S$) is an embedding,
 $\bar S$-close to $\bar\varphi$,
 that proves the lemma.
 \end{proof}

\begin {figure}
\includegraphics {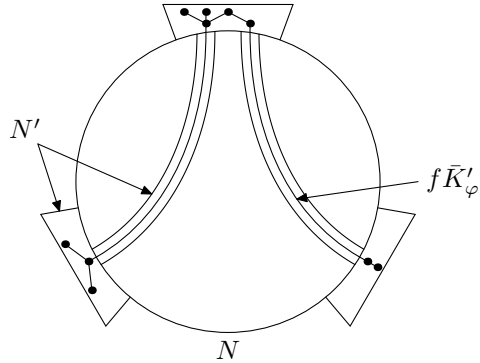}
\caption {Construction of the $S$-approximation}
\label {t-fig5}
\end {figure}

\begin {figure}
\includegraphics {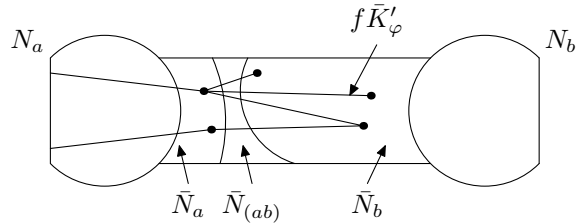}
\caption {Construction of the handle decomposition}
\label {t-fig6}
\end {figure}

 The same idea is used in the proof of Lemma~2.2.A,V.
 We take a general position map
 $f:\bar K'_\varphi\to N$, $S$-close to $\bar\varphi$,
 and construct its {\it semi-derivative}
 $\bar f':\bar K'_\varphi\to N'$, $S$-close
 to $\bar\varphi'$ (Fig.~7).
 Then we prove that if $f$ is an embedding, then
 $\bar f'$ is also an embedding (Fig.~8).

 \begin{definition*}[Definition~of~$\bar f'$, see Fig.~7,
 where this construction is applied to the map $\varphi$ from Fig.~4]
 Let $K$ be a graph without vertices of degree 0.
 Let $\varphi:K\to G\subset N$ be a nondegenerate
 simplicial map without transversal self-intersections.
 Let $f:K\to N$ be an $S$-approximation of $\varphi$.
 Then {\it the semi-derivative} $S'$-approximation
 $\bar f':\bar K'_\varphi\to N'$ is constructed as follows.
 For each edge $a\subset G$ fix a homeomorphism $h_a:N_a\to N'_{a'}$
 such that
 for each edge $b$ adjacent to $a$ we have
 $h_a(N_a\cap N_{(b)})\subset N'_{(a'b')}$.
 Define $\bar f'$ on each $\varphi$-component
 $\bar\alpha'\subset\bar K'_\varphi$ by the formula
 $\bar f'\left|_{\bar\alpha'}\right.=
 h_{\varphi \alpha} f\left|_{\alpha}\right.$
 Now let us define $\bar f'$ on each edge $xy\subset \bar K_\varphi'$
 joining two distinct $\varphi$-components $\bar X'$ and $\bar Y'$.
 Take an edge $a\subset \bar X'$ containing $x$.
 Identify $\bar X'$ and $X$ (see Definition of $\bar \varphi$ and
 $\bar \varphi'$).
 Then $a$ is also an edge of $K$ and $x$ is also a vertex of $K$.
 Denote by $\bar x$ the arc $a\cap f^{-1}N_{\varphi x}$.
 Define the arc $\bar y$ analogously.
 Decompose the edge $xy$ into three segments $xx_1$, $x_1y_1$ and $y_1y$.
 Let $\bar f'$ homeomorphically map $xx_1$ onto $h_{\varphi X} f\bar y$,
 $y_1y$ onto $h_{\varphi Y} f\bar x$, and
 $x_1y_1$ onto the rectilinear segment in
 $N'_{(\varphi X\,\varphi Y)}$ joining
 the  points $\bar f'x_1$ and $\bar f'y_1$.
 Thus the map $\bar f':\bar  K_\varphi'\to N'$
 is constructed.
 \end{definition*}

 \begin {figure}
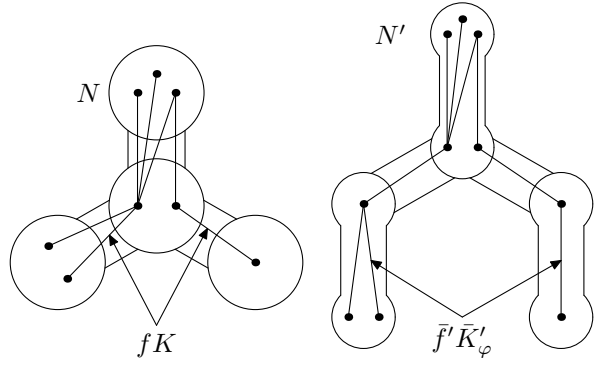

 \begin{tabular}{cc}
\includegraphics {skopenkov2.701} & \includegraphics {skopenkov2.702}
\end{tabular}
\caption {Semiderivative of an $S$-approximation}
\label {t-fig7}
\end {figure}

 Note that if $f$ is an embedding then
 there is a simpler alternative construction of $\bar f'$,
 in some sense reverse to the construction from the proof of Lemma~2.1.
 But this alternative construction is useless in the proof of Lemma~2.2.V,
 so we do not present it in the paper.
 We prove 2.2.A,V only in case
 when the derivative $N'$ is well-defined,
 i. e. $K$ does not contain pairs
 of arcs $i,j$ such that $\varphi i\cap\varphi j$ is transversal.
 This is sufficient for the proof of Theorem~1.3 and~1.5.
 In general case the proof is completely analogous, but one should
 use the definition of derivative $D(N)$ from \cite{6}.

\begin {figure}
\includegraphics {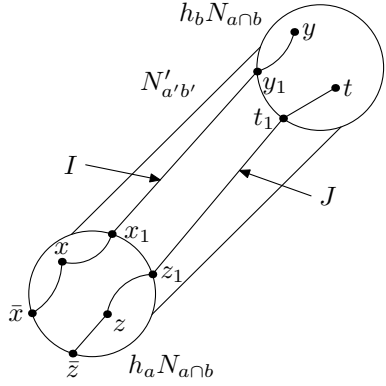}
\caption {Counting the number of crossings}
\label {t-fig8}
\end {figure}

 \begin{proof}[Proof of 2.2.A]
 By Contracting~Edge~Proposition~2.5.K
 we may assume that $\varphi$ is nondegenerate.
 Take an embedding $f:K\to N$, $S$-close to $\varphi$.
 Then it suffices to show that the map $\bar f'$
 (see Definition~of~$\bar f'$) is an embedding.

 Consider a pair of distinct edges $xy,zt$ of $K'_\varphi$.
 Denote the set $\bar f'(xy)\cap \bar f'(zt)$ by $i$.
 It suffices to show that $i=\bar f'(xy\cap zt)$.
 Denote by $a'=\bar\varphi'x$, $b'=\bar\varphi'y$,
 $c'=\bar\varphi'z$ and $d'=\bar\varphi't$.
 Without loss of generality we have the following 3 cases.

 1) $a'$, $b'$, $c'$ and $d'$ are pairwise distinct.
 Since $\bar f'$ is an $S'$-approximation, it follows that
 $\bar f'xy\subset N'_{a'b'}$ and
 $\bar f'zt\subset N'_{c'd'}$, hence $i=\emptyset$.

 2) ($a'=c'$ and $b'\ne d'$) or ($a'=b'=c'=d'$).
 Then $i\subset N'_{a'}$, hence $i=h_a(f\bar x\cap f\bar z)$
 (see the definition of $h_a$ and $\bar x$ in Definition~of~$\bar f'$,
 define $\bar z$ analogously to $\bar x$.)
 If $y\ne t$, then $\bar x$ and $\bar z$ are disjoint,
 so $f\bar x\cap f\bar z=\emptyset$ and $i=\emptyset$.
 If $y=t$, then  $i=h_a(fy)=\bar f'(xy\cap zt)$.

 3) $a'=c'$, $b'=d'$ and $a'\ne b'$.
 In this case both $xy$ and $zt$ join the vertices
 of distinct $\varphi$-components.
 Let us prove that $xy$ and $zt$ are disjoint.
 For example, assume that $y=t$.
 Then all the vertices $x$, $y$, $z$ and $t$ of $\bar K'_\varphi$
 correspond to the same vertex of $K$ denoted by $w$.
 Denote by $X$ and $Z$ the $\varphi$-components of
 $\varphi^{-1}a=\varphi^{-1}c$ such that $x\in \bar X'$ and $z\in\bar Z'$.
 So the $\varphi$-components $X$ and $Z$ have a common point $w$,
 hence $X=Z$.
 So $x,z\in \bar X'=\bar Z'$ correspond to the same vertex $w$,
 hence $x=z$.
 We obtain $y=t$ and $x=z$, then by the construction of $\bar K'_\varphi$
 we get $xy=zt$, that contradicts to the choice of these edges.
 So $xy$ and $zt$ are disjoint.

 Let us show that in case (3) $|i|=0 \pmod2$.
 Omit $\bar f'$ from the notation of $\bar f'$-images.
 Note that the homeomorphism $h_a\circ h_b^{-1}$
 maps $y_1y$ and $t_1t$ onto $\bar x$ and $\bar z$
 respectively (Fig.~8).
 First this implies that $|i|=|I\cap J|$,
 where $I=\bar x\cup xy_1$ and $J=\bar z\cup zt_1$.
 Secondly this implies that
 the two pairs of points $\partial I$ and
 $\partial J$ are not linked in $\partial(h_aN_{a\cap b}\cup N'_{(a'b')})$.
 Since $I,J\subset h_a N_{a\cap b}\cup N_{(a'b')}$,
 it follows that $|i|=|I\cap J|=0\pmod2$.
 So it remains to prove that $|I\cap J|\le1$, then $I\cap J=\emptyset$.
 This follows from
 $$
 \bar x\cap \bar z=h_a(f\bar x\cap f\bar z)=\emptyset\qquad
 xx_1\cap zz_1=h_a(f\bar y\cap f\bar t)=\emptyset
 \qquad\text{and}\qquad|x_1y_1\cap z_1t_1|\le1,
 $$
 because
 $x_1y_1$ and $z_1t_1$ are rectilinear segments in
 $N_{(a'b')}$.
 This completes the proof of the lemma.
 \end{proof}

 \begin{proof}[Proof of 2.2.V]
 By Contracting Edge Proposition~2.5.V
 it suffices to prove that
 if $f:\bar K'_\varphi\to N$ is a $\mod2$-embedding, $S$-close to $\varphi$,
 then the semi-derivative $\bar f'$ is also a $\mod2$-embedding.

 Take a pair of disjoint edges $xy,zt$ of $\bar\varphi'$
 and consider the three cases from the proof of Lemma~2.2.A.
 Case 1) is trivial.
 In case 2) we get $f(xy)\cap f(zt)\subset N_a$, hence
 $|i|=|h_a(f\bar x\cap f\bar z)|=|h_a(f(xy)\cap f(zt))|
 =|f(xy)\cap f(zt)|=0\pmod2$.
 In the proof of 2.2.A it is shown that in case 3)
 $|i|=0\pmod2$, thus Lemma~2.2.V is proved.
\end{proof}

 \section{The van Kampen obstruction}

 {\it The van Kampen obstruction} was
 invented by van Kampen in
 studies of embeddability of polyhedra in $\R^{2n}$ \cite{2, 3, 4,
 7, 8}.
 Let us give the definition of the van Kampen obstruction
 to approximability by embeddings of simplicial paths.
 Our construction is more visual than that in the problem of embeddability.
 Let $\varphi:I\to \R^2$ be a simplicial path (in Fig.~9 the constuction
 below is applied to the path from Fig.~1).
 Denote by $x_1,\dots,x_k$ the vertices of $I$ in the order along $I$,
 and denote the edge $x_ix_{i+1}$ by $i$.
 Let $I^*=\bigcup\limits_{i<j-1}i\times j$
 be {\it the deleted product} of $I$.
 Paint red the edges $x_i\times j$, $j\times x_i$, and the cells $i\times j$
 of $I^*$ such that $\varphi x_i\cap\varphi j=\emptyset$,
 $\varphi i\cap\varphi j=\emptyset$, and denote by $I^{*\varphi}$
 the red set.
 Take a general position map $f:I\to\R^2$, sufficently close to $\varphi$.
 To each cell $i\times j$ of ''the table'' $I^*$ put
 the number $v_f(i\times j)=|fi\cap fj|\pmod2$.
 Decompose $I^*$ along the red edges,
 let $C_1,C_2,\dots,C_n$ be all the obtained
 components such that
 $\partial C_k\cap\partial I^*\subset I^{*\varphi}$.
 Denote  $v_f(C_k)=\sum\limits_{i\times j\subset C_k}
 v_f(i\times j)$.
 {\it The van Kampen obstruction} (with $\Z_2$-coefficients) for
 approximability by embeddings is the vector
 $v(\varphi)=(\,v_f(C_1),v_f(C_2),\dots ,v_f(C_n)\,)$.

 \begin {figure}
\includegraphics {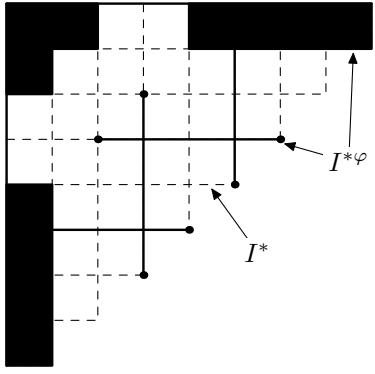}
\caption {The Van Kampen obstruction}
\label {t-fig9}
\end {figure}

 It can be shown easily that $v(\varphi)$ does not depend on the choice
 of $f$ \cite{8},
 thus $v(\varphi)=0$ is a necessary condition for approximability by
 embeddings.
 It is easy to check that $v(\varphi)\ne0$ for a PL path $\varphi:I\to\R^2$
 containing a transversal self-intersection.
 Thus Corollary~1.4.V follows from 1.3, 2.2.V and 3.1.

 \begin{proposition}
 The van Kampen obstruction $v(\varphi)=0$ if and only if
 there is an $S$-close
 to $\varphi$ general position $\mod2$-embedding.
 \end{proposition}

 \begin{proof}[Proof of 3.1]
 The inverse implication of the proposition is obvious.
 The proof of the direct implication follows the idea of [4].
 We are going to use
 the cohomogical formulation of the van Kampen obstruction
 (see the paragraph before Proposition~3.2 below for details).
 Let $f:K\to N$ be any general position $S$-approximation of $\varphi$.
 The 'Reidemeister move' shown in Fig.~10.a adds to $v_f$ the coboundary
 $\delta[x\times a]$ of the elementary cochain from $B^2(\tilde K)$.
 Since $v(\varphi)=0$, it follows that using some such 'moves' we
 can obtain a map $f:K\to N$ such that $v_f=0$.
 Then $f$ is the required $\mod2$-embedding, because $v_f=0$ yields that
 $|fa\cap fb|=0\pmod2$ for any pair of
 disjoint edges $a,b$ of $K$.
\end{proof}

 \begin {figure}
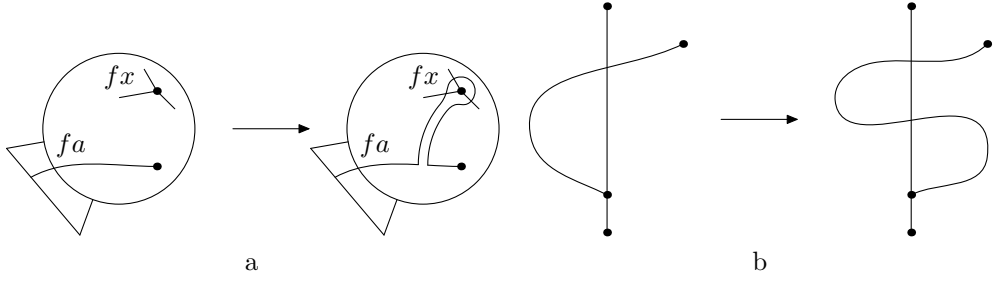

 \begin{tabular}{cc}
\includegraphics {skopenkov2.101} & \includegraphics {skopenkov2.102}\\
a & b
\end{tabular}
\caption {''The Reidemeister move''}
\label {t-fig10}
\end {figure}

 Now we are going to prove that the conditions 1.4.V and 1.4.D
 are equivalent (Proposition~3.2).
 We are going to replace $\Z_2$-coefficients in the van Kampen obstruction
 by $\Z$-coefficients, so Proposition~3.2 implies only that
 $1.4.D\implies 1.4.V$,
 but this is sufficient for the proof of Corollary~1.4.
 We prove Proposition~3.2 in the most general situation, so we need
 some more definitions.

 Let $K$ be an $n$-polyhedron with a fixed triangulation.
 Let $\varphi:K\to G\subset\R^{2n}$ be a simplicial map.
 Denote by $\sigma$ and $\tau$ any $n$-dimensional simplices
 of this triangulation of $K$.
 By {\it the deleted product} of $K$
 we mean the set
 $\tilde K=\bigcup\{\,\sigma\times\tau:\sigma\cap\tau=\emptyset\,\}$.
 Fix the natural orientation of each cell $\sigma\times\tau$
 (a positive basis of $\sigma\times\tau$ consists of
 the vectors $e_1,\dots,e_{2n}$, where $e_1,\dots,e_n$ form
 a positive basis of $\sigma$ and $e_{n+1},\dots,e_{2n}$ form
 a positive basis of $\tau$).
 Let $K^*=\tilde K/\Z_2$ be the factor under antipodal $\Z_2$-action.
 Let $\tilde K^{\varphi}\subset \tilde K$ be the subset
 $\tilde K^{\varphi}=
   \{\,\sigma\times\tau:\varphi\sigma\cap\varphi\tau=\emptyset\,\}$
 and let $K^{*\varphi}=\tilde K^{\varphi}/\Z_2$.
 For a general position map $f:K\to\R^{2n}$ close to $\varphi$
 define a cochain
 $v_f\in C^n(K^*,K^{*\varphi};\Z)$ by the formula
 $v_f(\sigma\times \tau)=f\sigma\cap f\tau$.
 This cochain is well-defined, because
 $f\sigma\cap f\tau=(-1)^n f\tau\cap f\sigma$ and
 our $\Z_2$-action maps $\sigma\times\tau$ to $(-1)^n \tau\times\sigma$.
 The class $v(\varphi)=[v_f]\in H^n(K^*,K^{*\varphi};\Z)$ of this cochain
 does not depend on the map $f$ and is called {\it the van Kampen obstruction}
 to approximability of $\varphi$ by embeddings.
 We say that the map $\varphi:K\to G\subset\R^{2n}$ satisfies
 {\it the deleted product property} if the map
 $\tilde\varphi:\tilde K^{\varphi}\to S^{2n-1}$ given by the formula
 $\tilde\varphi(x,y)=\frac{\varphi x-\varphi y}{\| \varphi x-\varphi y \|}$
 extends to an equivariant map $\tilde K\to S^{2n-1}$.
 Evidently, this definition of the deleted product property
 is equivalent to 1.4.D for $n=1$ and $K\cong I$.

 \begin{proposition}
 A PL map $\varphi:K^n\to \R^{2n}$ satisfies the deleted product property
 if and only if the van Kampen obstruction (with $\Z$-coefficients)
 $v(\varphi)=0$.
 \end{proposition}

 \begin{proof}[Proof of 3.2]
 We are going to show that the van Kampen obstruction
 is a complete obstruction to an equivariant extension of
 $\tilde\varphi:\tilde K^{\varphi}\to S^{2n-1}$
 to a map $\tilde K\to S^{2n-1}$.

 Take a general position map $f:K\to\R^{2n}$ close to $\varphi$
 and define the equivariant map
 $\tilde f:\tilde K^{\varphi}\cup\sk^{2n-1}\tilde K\to S^{2n-1}$ by
 the formula  $\tilde f(x,y)=\frac{fx-fy}{|fx-fy|}$.
 By general position it follows that $\tilde f$ is well-defined.
 Since $f$ is close to $\varphi$, it follows that
 $\tilde f\left |_{\tilde K^{\varphi}}\right.$ is homotopic to
 $\tilde \varphi$.
 Evidently, then $\tilde\varphi$ extends to an equivariant
 map $\tilde K\to S^{2n-1}$ if and only if
 $\tilde f\left |_{\tilde K^{\varphi}}\right.$ extends to an equivariant
 map $\tilde K\to S^{2n-1}$.

 Consider a cell $\sigma\times\tau\subset\tilde K-\tilde K^{\varphi}$,
 where $\sigma,\tau\subset K$ are $n$-dimensional cells.
 The map $\tilde f$ extends to $\sigma\times\tau$ if and only if
 $\operatorname{deg} \tilde f\left|{\partial(\sigma\times\tau)}\right. =0$.
 If $\tilde f$ extends to $\sigma\times\tau$ then it extends
 also to $\tau\times\sigma$ in equivariant way,
 because $\tilde f$ is equivariant.
 One can see that
 $\operatorname{deg} \tilde f\left|{\partial(\sigma\times\tau)}\right.
 =f\sigma\cap f\tau=v_f(\sigma\times \tau)$.
 So the map $\tilde f$ extends to an equivariant map
 $\tilde K\to S^{n-1}$ if and only if $v_f=0$.

 Now let $g:\tilde K^{\varphi}\cup\sk^{2n-1}\tilde K\to S^{2n-1}$
 be an equivariant map such that $gx=\tilde fx$ for each
 $x\in\tilde K^{\varphi}\cup\sk^{2n-2}\tilde K$.
 Define the cochain $v_g\in C^{2n}(K^*,K^{*\varphi};\Z)$
 by the formula
 $v_g(\sigma)=\deg g\left|{\partial\sigma}\right.$
 for each $2n$-dimensional cell $\sigma$.
 Let $\sigma\subset\tilde K-\tilde K^{\varphi}$ be
 a cell of dimension $2n-1$.
 Take a disjoint union $\sigma\sqcup\sigma'$ of two copies of $\sigma$
 and paste $\sigma$ to $\sigma'$ by $\partial\sigma=\partial\sigma'$.
 Let $d_{\sigma}$ be a map of the obtained $(2n-1)$-sphere to $S^{2n-1}$
 given by the formula $d_{\sigma}x=fx$ for $x\in\sigma$ and
 $d_{\sigma}x=gx$ for $x\in\sigma'$.
 Define the cochain $v_{fg}\in C^{2n-1}(K^*,K^{*\varphi};\Z)$
 by the formula $v_{fg}(\sigma)=\operatorname{deg} d_{\sigma}$
 (we fix the orientation of $\sigma\cup\sigma'$ restricted
 to the positive orientation of $\sigma'$).
 Then, clearly $v_g-v_f=\delta v_{fg}$.

 The obtained formula implies that the cohomological class $[v_g]$
 does not depend on the choice of an equivariant map
 $g:\tilde K^{\varphi}\cup\sk^{2n-1}\tilde K\to S^{2n-1}$ and
 coincides with the van Kampen obstruction $v(\varphi)$
 (with $\Z$-coefficients).
 This proves that the condition $v(\varphi)=0$ in the proposition
 is necessary.
 The obtained formula and the construction of $v_{fg}$ above shows
 that if $v(\varphi)=0$ then $v_g=0$ for some
 $g:\tilde K^{\varphi}\cup\sk^{2n-1}\tilde K\to S^{2n-1}$,
 hence $\tilde f\left|_{\tilde K^{\varphi}}\right.$ extends
 to an eguivariant map $\tilde K\to S^{2n-1}$.
 So the proposition is proved.
 \end{proof}

 \begin{example} (cf. \cite{15, 1})
 There exists a pair of PL paths $\varphi:I\to \R^2$,
 $\psi:I\to \R^2$
 (Fig.~11, where a pair of paths $f$, $g$, close to them, is shown),
 not approximable by link maps (i. e., maps with disjoint images)
 and such that:

 V) {\it The van Kampen obstruction} $v(\varphi,\psi)=0$.

 D) The map
 $\Phi:\{\,(x,y)\in I\times I\,|\,\varphi x\ne\psi y\,\}\to S^1$
 given by $\Phi(x,y)=\frac {\varphi x-\psi y}
 {\|\varphi x-\psi y\|}$ homotopically extends to a map $I\times I\to S^1$.

 I) The pair $\varphi'$, $\psi'$ is approximable by link maps.
 \end{example}

 \begin {figure}
\includegraphics {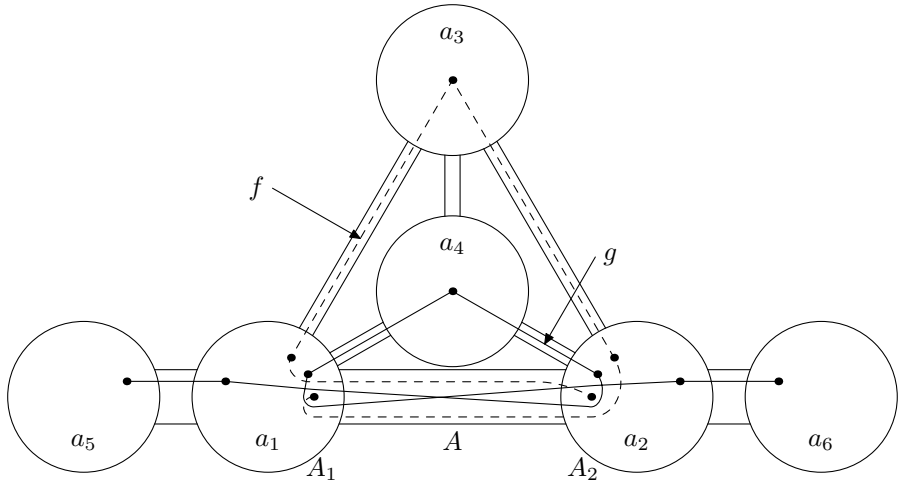}
\caption {A pair of maps not approximable by link maps}
\label {t-fig11}
\end {figure}

 \begin{proof}[Proof of 3.3]
 Let $K,L\cong I$ are the graphs with the vertices $k_1,\dots,k_5$
 and $l_1,\dots,l_7$, let $G$ be the graph with vertices $a_1,\dots,a_6$
 and edges $a_1a_2$, $a_1a_3$, $a_1a_4$, $a_1a_5$, $a_2a_3$,
 $a_2a_4$ and $a_2a_6$.
 The required simplicial maps $\varphi,\psi$ are given by the formulae
 $\varphi k_1=a_1$, $\varphi k_2=a_2$, $\varphi k_3=a_3$, $\varphi k_4=a_1$,
 $\varphi k_5=a_2$ and $\psi l_1=a_5$, $\psi l_2=a_1$, $\psi l_3=a_2$,
 $\psi l_4=a_4$, $\psi l_5=a_1$, $\psi l_6=a_2$, $\psi l_7=a_6$.
 Consider the pair of $S$-approximations $f$ and $g$ of $\varphi$ and $\psi$
 respectively shown in Fig.~11.
 One can see that $|fi\cap gj|=0\pmod2$ for any pair of edges
 $i\subset K, j\subset L$.
 This implies both 3.3.V and 3.3.D (it is shown analogously
 to the proof of 1.4, see also Proposition~3.1).
 The proof of 3.3.I is a direct calculation.
 Let us prove that the pair $\varphi,\psi$ is not approximable by
 link maps.  Assume the converse.
 Let $K_{13},K_{35}\subset K$ and $L_{14},L_{47}\subset L$
 be the arcs between the points $k_1$ and $k_3$, $k_3$ and $k_5$,
 $l_1$ and $l_4$, $l_4$ and $l_7$ respectively.
 Take a small neighbourhood of $\varphi K\cup \psi L$
 in the plane and fix its handle decomposition $S$.
 Denote by $A_1$, $A_2$ and $A$ the discs of $S$
 corresponding to the vertices $a_1$, $a_2$ and to the edge $a_1a_2$
 respectively.
 By the analogue of the Minc Proposition
 (see the paragraph after Definition of an~$S$-approximation in \S2)
 there are two $S$-approximations $f,g$ of
 $\varphi$ and $\psi$ respectively, having disjoint images.
 Since $fK_{13}\cap gL=\emptyset$, it follows that
 the pairs of points $gL_{14}\cap\partial (A_1\cup A)$
 and $gL_{47}\cap\partial (A_1\cup A)$ are not linked in
 $\partial(A_1\cup A)$.
 Analogously, $gL_{14}\cap\partial A_2$ and $gL_{47}\cap\partial A_2$
 are not linked in $\partial A_2$.
 So $gL_{14}\cap\partial(A_1\cup A_2\cup A)$ and
 $gL_{47}\cap\partial(A_1\cup A_2\cup A)$
 are not linked in $\partial(A_1\cup A_2\cup A)$.
 Then $g$ cannot be an $S$-approximation of $\psi$.
 This contradiction proves that $\varphi$ and $\psi$
 are not approximable by link maps.
\end{proof}

 \section{Variations}

 The following Simple-minded Criterion~4.1
 for approximability by embeddings gives an algorithm
 of checking whether a given nondegenerate map is approximable by
 embeddings (another algorithm is given in \cite{13}).

 \begin{proposition}[Simple-minded Criterion] Let $\varphi:K\to G\subset\R^2$ be
 a nondegenerate simplicial map of a graph $K$,
 i. e. for each edge $a$ of $K$ the image $\varphi a$ is not
 a vertex.
 Replace each edge $a\subset G$ by $i$ close multiple edges in $\R^2$,
 if $\varphi^{-1}a$ consists of $i$ edges.
 Denote by $\bar G\subset \R^2$ the obtained graph and by
 $\pi:\bar G\to G$ the evidently defined projection.
 The map $\varphi$ is approximable by embeddings
 if and only if there exists
 an onto map $\bar\varphi:K\to\bar G$ without transversal self-intersections
 and such that $\pi\circ \bar\varphi=\varphi$.
 \end{proposition}

 The proof is trivial (we do not present the details since we do not use
 this criterion).
 There exists a purely combinatorial proof of Theorem~1.3,
 based on Criterion~4.1.
 Criterion~4.1 and all the other our previous results remain
 true, if we replace $\R^2$ by an arbitrary orientable 2-manifold $N$.

 There exists an infinite number of PL maps
 $\varphi:T\to T\subset\R^2$, where $T$ is letter ''T''
 ({\it a simple triod}), not approximable by embeddings and
 such that $\varphi'$ and any simplicial restriction of $\varphi$ are
 approximable by embeddings (for the only embedding $T'\to\R^2$).
 So there are no criteria like 1.3.I,S for $K\not\cong I,S^1$.

 In the rest of the paper we prove Theorem~1.5,
 generalizing both 1.3 and 1.4.
 For the proof we need the following
 Lemma~4.2.T, Lemma~4.5 and Lemma~4.6,
 analogous to Lemmas~2.3, 2.1 and 2.2 respectively.

 To state Lemma.~4.2 we need the following definitions.
 We shall say that $\varphi$
 {\it contains a simple triod}, if there is a triod $T\subset K$
 with the edges $t_1, t_2, t_3$ such that
 the arcs $\varphi t_1,\varphi t_2,\varphi t_3$
 have a unique common point.
 We shall say that $\varphi$
 {\it identifies triods}, if it contains two disjoint simple
 triods $T_1, T_2\subset K$ such that $\varphi T_1=\varphi T_2$.
 Note that $v(\varphi)\ne0$ for a map $\varphi$ identifying
 triods.
 We shall say that an onto map $\varphi:K\to G$ is {\it a standard winding},
 if both $K$ and $G$ are
 homeomorphic to disjoint unions of circles (may be,
 $K=G=\emptyset$) and $\varphi\left|_S\right.$ is
 a standard winding of a nonzero degree for each circle
 $S\subset K$.
 Denote by
 $\Sigma(\varphi)=\{\, x\in K \,:\, |\varphi^{-1}\varphi x|\ge2 \,\}$
 {\it the singular set} of the map $\varphi$.

 \begin{lemma} (cf. Lemma~2.3)
 Let $\varphi:K\to G$ be a simplicial map of a graph $K$ with $k$
 vertices.

 T) Suppose that for each $i=0,\dots,k$ the derivative
 $\varphi^{(i)}$ does not contain simple triods;
 then $\varphi^{(k)}$ is a standard winding.

 I) Suppose that for each $i=0,\dots,k$ the derivative
 $\varphi^{(i)}$ does not identify triods;
 then $\Sigma(\varphi^{(k)})$ is a disjoint union of circles and
 $\varphi^{(k)}\left|_{\Sigma(\varphi^{(k)})}\right.$
 is a standard winding.
 \end{lemma}

 \begin{proof}[Proof of 4.2.T]
 We are going to use the notation from the proof of Lemma~2.3.
 Let us show that $|K|\ge |K_\varphi'|$ for a simplicial map
 $\varphi:K\to G$ containing no simple triods.
 We also prove that
 $|K|=|K_\varphi'|$ only if $K$ is homeomorphic
 to a disjoint (maybe empty) union of circles
 and $\varphi$ is ultra-nondegenerate.
 Then it suffices to prove the lemma for this latter specific case.
 Indeed, since $\varphi$ does not contain simple triods,
 it follows that each vertex of the graph $K$ belongs to at most
 two $\varphi$-components of the graph.
 On the other hand, each $\varphi$-component contains an edge,
 and hence it contains at least two vertices.
 This yields that the number of vertices of $K$
 is greater or equals to
 the number of $\varphi$-components, i. e. $|K|\ge |K_\varphi'|$.
 We have the equation here if and only if (1) each vertex of $K$ belongs to
 two $\varphi$-components and (2) each $\varphi$-component
 contains exactly two vertices.
 The condition (2) means that $\varphi$ is ultra-nondegenerate.
 But for an ultra-nondegenerate map the condition (1)
 yields that the degree of each
 vertex of $K$ is 2, so $K$ is a disjoint union of circles.
 Now note that for any two components $A,B\subset K$
 we have $\varphi'A'\cap \varphi'B'=(\varphi A\cap \varphi B)'$.
 Since $A,B\cong S^1$ and $\varphi$ is ultra-nondegenerate,
 it follows that
 $\varphi A\cap \varphi B$ is always either empty or a circle or
 a disjoint union of arcs and points.
 Moreover, $\varphi A\cap\varphi B$ is a circle if and only if
 $\varphi A=\varphi B$.
 If $\varphi A\cap\varphi B$ is either empty or a disjoint union of arcs
 and points,
 then $\varphi^{(k)}A$ and $\varphi^{(k)}B$ are disjoint.
 So the images $\varphi^{(k)}A$ and
 $\varphi^{(k)}B$ either are disjoint or coincide.
 By Lemma~2.3 this yields that $\varphi^{(k)}$ is
 a standard winding.
 \end{proof}

 We do not use Lemma~4.2.I
 and prove it after the proof of Theorem~1.5.
 Lemma~4.2.I may be helpful in the proof of Conjecture~1.6.

 In order to prove Theorem~1.5 we need the following extension of the
 techniques from \S2.
 As in \S2, we fix an orientable thickening $N$ of the graph $G\cong S^1$.
 We also assume that some orientable thickening $M$
 of the graph $K$ is fixed.
 Note that a thickening of a graph is uniquely defined by
 {\it a local ordering} of edges around each vertex \cite{6}.
 We assume that $K$ may contain loops and multiple edges.
 In this case by {\it a simplicial map} $\varphi:K\to G$
 we mean a continuous map that is linear on each edge of $K$ and such that
 $\varphi x$ is a vertex of $G$ for each vertex $x\in K$.
 We assume that the handle decompositions
 $M=\bigcup M_x\cup\bigcup M_{(a)}$ and
 $N=\bigcup N_x\cup\bigcup N_{(a)}$
 are fixed (in both formulae the first union is over all vertices $x$
 and the second --- over all edges $a$).
 By $M_\alpha$ and $N_\beta$ we denote the restriction
 of the thickenings $M$ and $N$
 to subgraphs $\alpha\subset K$ and $\beta\subset G$ respectively.
 By {\it an $S$-approximation} of $\varphi$
 (or {\it $S$-close to $\varphi$} map) we mean
 a general position map $f:M\to N$ such that
 for any vertex $x\in K$ or edge $x\subset K$ we have
 $f M_x\subset N_{\varphi x}$
 and for any edge $x\subset K$ with nondegenerate $\varphi x$
 the set $M_x\cap f^{-1}N_{(\varphi x)}$
 is connected (cf. Definition of~$S$-approximation in \S2).
 If there is an $S$-close to $\varphi$ embedding $M\to N$ then
 we shall say that $\varphi$ {\it is approximable by embeddings $M\to N$}.
 By {\it a $\mod2$-embedding} we mean
 a general position map $f:M\to N$ such that $f\left|_{M_x}\right.$
 is an embedding for each vertex $x\in K$, $f\left|_{M_{(a)}}\right.$ is an immersion for each edge $a\subset K$ and
 $|fa\cap fb-f(a\cap b)|=0\pmod2$ for any pair of distinct edges
 $a,b\subset K\subset M$.
 The last notion appears in the following generalization of
 Proposition~3.1.

 \begin{lemma}
 Let $K$ be a graph such that the degree of each vertex of $K$ is at most 3.
 Let $\varphi:K\to G\subset\R^2$ be a simplicial map such that $v(\varphi)=0$.
 Then there exist a  $\mod2$-embedding $M\to N$, $S$-close to $\varphi$,
 for some thickenings $M$ and $N\subset\R^2$ of the graphs $K$ and $G$
 respectively.
 \end{lemma}

 \begin{proof}[Proof of 4.3] \cite{10}
 Let $N$ be a regular neighbourhood of $G$ in $\R^2$.
 Let $f:K\to N$ be the map given by Proposition~3.1.
  Since the degree of each vertex of $K$ is at most 3, it follows that
 we can remove intersections of adjacent edges, using the moves shown in
 Fig.~10.b.
 The obtained general position map
 $K\to N$ uniquely defines a thickening $M$ of $K$
 and extends to the required $\mod2$-embedding, $S$-close to $\varphi$.
\end{proof}

 In \cite{10} the move shown in Fig.~10.b is assumed to work
 for vertices of any degree, that is not right.
 The degree restriction in Theorem~1.5
 is used only in this step of the proof.

 Now we are going to construct {\it the derivative $M'_\varphi$}
 of the thickening $M$.
 This derivative is well-defined only under the following
 conditions on $\varphi$, $M$ and $N$.
 We shall say that $\varphi:K\to G$
 is {\it locally approximable by embeddings},
 if for each vertex $x\in G$ there exists an $S$-approximation $f:M\to N$ of
 $\varphi$ such that $f\left|_{f^{-1}N_x}\right.$ is an embedding.
 The following lemma asserts that
 the thickenings $M$ and $N$ given by Lemma~4.3
 satisfy these conditions.

 \begin{lemma}[Local Approximation Lemma]
 If there is a $\mod2$-embedding $M\to N$ which is $S$-close to the map
 $\varphi:K\to G$,
 then $\varphi$ is locally approximable by embeddings $M\to N$.
 \end{lemma}

 \begin{proof}[Proof of 4.4]
Let $f:M\to N$ be an $S$-close to $\varphi$ $\mod 2$-embedding and
let $x$ be a vertex of the graph $G$.
Modify $f$ in a small neighbourhood of $\partial N_x$ to obtain a map $f$
such that for each edge $a\ni x$ the set $f^{-1}(N_{(a)}\cap N_x)$
is a single point $P_a$.
Attach a ring $R$ to the disc $N_x$ along the circle $\partial N_x$.
Let $Q$ be the 2-polyhedron obtained from
$M_{\varphi^{-1}x}\cup(K\cap f^{-1}N_x)$ by identifying the points
$f^{-1}P_a$ for each $a\ni x$.
Identify each point $f^{-1}P_a\in Q$ with $P_a$.
Attach the ring $R$ to $Q$ by the inclusions $P_a\subset R$.
Let $g:Q\cup R\to N_x\cup R$ be the map given by the formulae
$gy=fy$ for $y\in Q$ and $gy=y$ for $y\in R$.
Clearly, $g$ is {\it a $\mod 2$-embedding}, i. e. there exists
a triangulation of $Q\cup R$ such that $|ga\cap gb|=0\pmod2$
for each pair of disjoint edges $a,b$ of this triangulation.
This yields that $Q\cup R$ contains neither Kuratowsky graph $K_5$
nor $K_{3,3}$. Clearly, $Q\cup R$ contains neither $S^2$ nor the cone over
$S^1\sqcup D^0$.
By the well-known 2-polyhedron planarity criterion $Q\cup R$ is planar.
So the embedding $R\subset N_x\cup R$ extends to an embedding
$h:Q\cup R\to N_x\cup R$.
Clearly, $h\left|_Q\right.:Q\to N_x$ can be modified to an embedding
$f^{-1}N_x\to N_x$, that extends to the required $S$-approximation
of $\varphi$.
So $\varphi$ is locally approximable by embeddings.
\end{proof}

 The next step in the construction of $M'_\varphi$
 is like Contracting~Edge~Proposition~2.5.
 We use a reduction to the case of a nondegenerate map $\varphi$
 and then define the semi-derivative thickening.

 \begin{definition*}[Definition of $\varphi^c$]
 First let us define the graph $K_\varphi^c$,
 which is the domain of $\varphi^c$.
 A {\it $0$-component} of $K$ is any connected component of $\varphi^{-1}x$
 for a vertex $x\in G$.
 The vertex set of the graph $K_\varphi^c$ is in 1-1 correspondence
 with the set of all $0$-components.
 Denote by $\alpha^c$ the vertex corresponding to a 0-component
 $\alpha\subset K$.
 The vertices $\alpha^c$ and $\beta^c$ are joined by an edge in
 $K_\varphi^c$ if and only
 if $K$ contains an edge with two ends belonging to $\alpha$ and $\beta$
 respectively.
 The map $\varphi^c:K_\varphi^c\to G$ is a simplicial map
 given by the formula $\varphi^c\alpha^c=\varphi\alpha$.
 \end{definition*}

 \begin{definition*}[Definition of $M^c$] Let the map $\varphi:K\to G$
 be locally approximable by embeddings.
 The thickening $M_\varphi^c$
 and its handle decomposition are defined as follows.
 For each $0$-component $\alpha\subset K$ choose a maximal tree
 $T\subset \alpha$.
 The thickening $M_\varphi^c$ is the restriction of $M$
 to the subgraph $(K-\bigcup \alpha)\cup\bigcup T$.
 The discs of the handle decomposition of $M_\varphi^c$ are defined as
 the subthickenings $M_T$ and the strips are defined as the strips of $M$
 not contained in the subthickenings $M_{\alpha}$ and $M_T$.
 \end{definition*}

 \begin{definition*}[Definition of $\bar M'$,
 cf. Definition of $\bar\varphi'$ in \S2]
 Let the map $\varphi:K\to G$
 be nondegenerate and locally approximable by embeddings.
 The thickening $\bar M_\varphi'$
 and its handle decomposition are defined as follows.
 Take a disjoint union of $M_\alpha$ for all
 $\varphi$-components $\alpha\subset K$.
 If $\varphi$-components $\alpha$ and $\beta$ have a common vertex $x$,
 then join $M_\alpha$ and $M_\beta$ by a strip attached to
 $M_\alpha$ and $M_\beta$
 along the arcs $M_x\cap M_{(a)}$ and $M_x\cap M_{(b)}$
 respectively, where $a\subset\alpha$ and $b\subset\beta$ are any edges
 containing $x$.
 The handle decomposition of $\bar M_\varphi'$ is obtained from
 those of $M_\alpha$ by adding the attached strips to the decomposition.
 \end{definition*}

 We omit $\varphi$ from the notation of $M^c_\varphi$,
 $\bar M'_\varphi$ and $M'_\varphi$.
 We denote {\it the derivative} of the thickening $M$
 by $M'=((\bar M^{c})')^c$.
 Note that $M'$ is a thickening of a graph which may be
 different from $K_\varphi'$ only in multiplicity of some edges.
 Now we are going to prove
 the following Lemma~4.5.C,D and Lemma~4.6.C,D,
 analogous to 2.5.A, 2.1, 2.5.V and 2.2.V respectively.

 \begin{lemma}
 Let $\varphi:K\to G$ be locally approximable by embeddings. Then:

 C) $\varphi^c$ is approximable by embeddings $M^c\to N$ if and only if
 $\varphi$ is approximable by embeddings $M\to N$;

 D) if $G\cong I$ or $G\cong S^1$ and
 $\bar \varphi'$ is approximable by embeddings $\bar M'\to N'$,
 then $\bar\varphi$ is approximable by embeddings $\bar M'\to N$.
 \end{lemma}

 \begin{proof}[Proof of 4.5.C]
 The sufficiency is obvious because $M^c\subset M$
 (without handle decompositions).
 Let us prove the necessity.
 Let $f:M^c\to N$
 be an embedding, $S$-close to $\varphi^c$.
 Identify $M^c$ and the subthickening $M_{\bar K}$,
 where $\bar K$ is the unionwhere $\bar K=(K-\bigcup\alpha)\cup\bigcup T$
(see the notation in Definition of~$M^c$).
 Let us add to $\bar K$ the edges $a$ from $\alpha-T$
 one by one and
 extend the map $f$ to the corresponding strips $M_{(a)}$
 in an arbitrary way.
 We assert that an embedding $M\to N$, $S$-close to $\varphi$,
 can be constructed in this way.
 Indeed, assume that algorithm does not work, i. e.
 at some step the obtained $S$-close to $\varphi$ embedding
 $M_{\bar K}\to N$ cannot be extended to $M_{(a)}$.
 Then the two arcs $f\partial M_{(a)}$  are contained in
 distinct connected components of $\Cl(N_{\varphi\alpha}-fM_{\bar K})$.
 Since the arcs $f\partial M_{(a)}$ belong to one connected
 component of $N_{\varphi\alpha}\cap fM_{\bar K}$, it follows that
 there are only two possibilities:
 1) $\bar K$ contains a cycle $\gamma$ such that
 $\varphi\gamma=\varphi\alpha$ and
 $M_{(a)}$ joins the two components of $\partial M_{\gamma}$
 (Fig.~12.a);
 2) $\bar K$ contains an arc $\gamma$ such that
 $\varphi\alpha\not\in\varphi\,\partial\gamma$ and
 $M_{(a)}$ joins the two components of
 $\partial M_{\gamma}-M_{\partial\gamma}$ (Fig.~12.b).
 Clearly, in both cases 1) and 2) the map $\varphi$ is not locally
 approximable by embeddings.
 This contradicts the assumption of the lemma,
 so the algorithm above gives us the required embedding $M\to N$,
 $S$-close to $\varphi$.
\end{proof}

 \begin {figure}
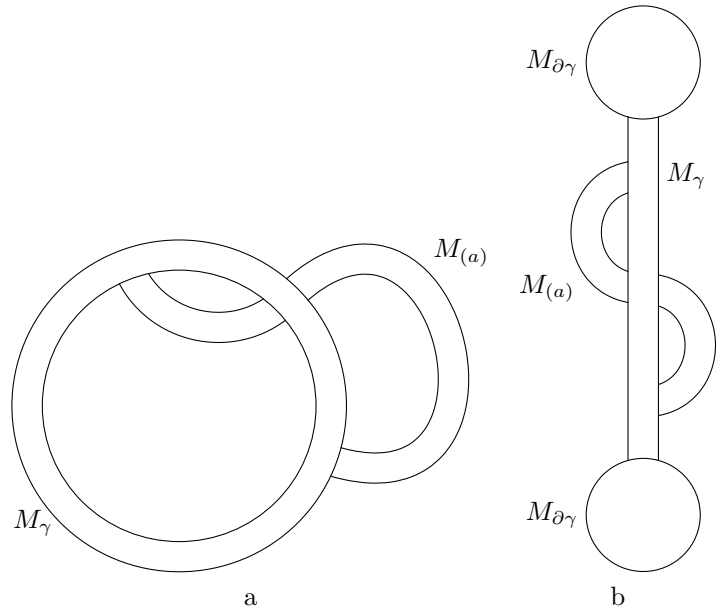

 \begin{tabular}{cc}
\includegraphics {skopenkov2.121} & \includegraphics {skopenkov2.122}\\
a & b
\end{tabular}
\caption {The thickenings}
\label {t-fig12}
\end {figure}

 \begin{proof}[Proof of Lemma.4.5.D]
 Note that $N'$ is well-defined because $G\cong S^1$
 and hence $\varphi i\cap\varphi j$ is not transveral
 for any pair of arcs $i,j\subset K$.
 Take an $S$-close to $\bar \varphi'$ embedding
 $\bar M'\to N'$.
 Let $f:\bar M'\to N$ be the composition of this embedding
 with the embedding $N'\to N$ (see Definition~of~$N'$).
 It remains to construct a new handle decomposition $\bar S$ of $N$
 such that $f$ is an $\bar S$-approximation of $\bar\varphi$.

  Denote by $m=\bar M'$.
 For each vertex $x\in G$ let $\bar N_x$ be a small neighbourhood in $N$
 of the set $f m_{\bar\varphi^{-1}x}\cup N_x$.
 Since $G\cong S^1$, it follows that $\bar\varphi^{-1}x$ is
 a disjoint union of arcs and hence $\bar N_x$ is a disjoint union of discs.
 For each edge $xy\subset G$ the set $f m_{\bar\varphi^{-1}xy} - \bar N_x$
 is also a disjoint union of discs, because
 for each edge $a$ such that $\bar\varphi a=xy$ the strip
 $f m_{(a)}$ joins $\bar N_y$ and some disc of $\bar N_x$.
 Hence $f m_{\bar\varphi^{-1}xy} - \bar N_x$ does not decompose $N$,
 and since $\bar N_x$ is a disjoint union of discs, it follows that
 $\bar N_y \cup(f m - \bar N_x)$ does not decompose $N_{(xy)}$.
 So if $\bar N_x$ is not connected, one can join
 any of its connected components $C\subset N_{(xy)}$
 with another connected component of $\bar N_x$ by a strip in $N_{(xy)}$ not intersecting $fm$
and $\bar N_y$.
 Clearly, $\bar N_x$ remains to be a disjoint union of discs
 after this operation.
 Let us add to each $\bar N_x$ such strips until
 all $\bar N_x$ become connected.
 Then the obtained discs $\bar N_x$ and the strips
 $\bar N_{(xy)}=\Cl(N_{(xy)}-\bar N_x-\bar N_y)$
 form the required handle decomposition $\bar S$.
\end{proof}

 \begin{lemma}
 If there is  a $\mod2$-embedding $M\to N$ $S$-close to $\varphi$,
 then there is a $\mod2$-embedding, $S$-close to
 \quad C) $\varphi^c$; \quad D) $\bar \varphi'$.
 \end{lemma}

 \begin{proof}[Proof of 4.6.C]
 Take an $S$-close to $\varphi$ $\mod2$-embedding $M\to N$.
 Denote by $f:K\to N$ the restriction of the embedding.
 Make the move from Contracting Edge Proposition~2.5.K (Fig.~3(b)) for
 each edge of the tree $T\subset K$ (see Definition of~$M^c$).
 Let $f:K_\varphi^c\to N$ be the restriction
 of the constructed map.
 By the construction the local ordering of the edges of $K_\varphi$
 in $M^c$ and the local ordering of their $f$-images coincide.
 So the map $f:K_\varphi^c\to N$ extends to the required
 $\mod2$-embedding $M^c\to N$.
 \end{proof}

 We prove Lemma 4.6.D only in case when $K$ does not contain
 pairs of arcs $i,j$ with $\varphi i\cap\varphi j$
 transversal (this assumption for $G\cong S^1$
 is satisfied automatically.) The lemma is proved analogously
 for an arbitrary graph $K$, if we use the definition of $N'$
 from \cite{6}.

 \begin{proof}[Proof of 4.6.D]
 Take an $S$-close to $\varphi$ $\mod2$-embedding $M\to N$.
 Denote by $f:K\to N$ the restriction of the embedding.
 Let $\bar f':\bar K_\varphi'\to N'$ be the map
 from Definition of~$\bar f'$ in~\S2.
 By the construction of Definition of~$\bar f'$ and Definition of $\bar M'$
 the local ordering of the edges of $\bar K_\varphi'$
 in $\bar M'$ and the local ordering of their $\bar f'$-images coincide.
 So the map $\bar f'$ extends to the required
 $\mod2$-embedding $\bar M'\to N'$.
\end{proof}

 \begin{proof}[Proof of 1.5]
 The necessity follows from the necessity of the condition
 $v(\varphi)=0$ for approximability by embeddings \cite{8},
 Example~1.1 \cite{11} and Lemma~2.2.A \cite{6}.
 Let us prove the sufficiency.
 Suppose that $v(\varphi)=0$.
 By Lemma~4.3 there exist a thickening $M$ of $K$
 and a $\mod2$-embedding $f:M\to N$, $S$-close to $\varphi$.
 Local Approximation Lemma~4.4 implies that $\varphi$ is locally approximable
 by embeddings.
 Note that the graphs $K_\varphi'$ and $((\bar K^c)')^c$ are the same
 modulo multiplicity of some edges,
 $\varphi'$ and $((\bar \varphi^c)')^c$ coincide modulo this difference.
 Change $\varphi'$ to this {\it quasi-derivative} $((\bar \varphi^c)')^c$.
 By Lemma~4.6.C,D it follows that for each natural $i$ the
 derivative thickening $M^{(i)}$ is well-defined
 and there exists an $S$-close to the map $\varphi^{(i)}$
 $\mod2$-embedding, and the map $\varphi^{(i)}$
 is locally approximable by embeddings.
 By Lemma~4.2.T {\it the derivative} $\varphi^{(k)}$  is either ''empty'' or
 the standard $d$-winding for some $d\ne 0$, hence
 {\it the quasi-derivative} $\varphi^{(k+1)}$
 is either ''empty'' or the standard $d$-winding with the same $d$.
 By the assumption of the theorem $d$ is either even or $\pm1$.
 Since for the standard winding of
 an even degree $d\ne0$ the van Kampen obstruction is nonzero
 and $\varphi^{(k+1)}$
 is approximable by $\mod2$-embeddings, it follows that $d$ is not even.
 Hence either $d=\pm1$ or $\varphi^{(k+1)}$
 has an empty domain. In both cases there exists an $S$-close to the map
 $\varphi^{(k+1)}$ embedding $K^{(k+1)}\to N^{(k+1)}$,
 where {\it the quasi-derivatives} $K^{(i)}$ are defined analogously
 to $\varphi^{(i)}$.
 Since $K^{(k+1)}=\emptyset$ or $K^{(k+1)}\cong S^1$,
 it follows that this embedding
 extends to an embedding $M^{(k+1)}\to N^{(k+1)}$.
 Applying Lemma~4.5.C,D $k+1$ times, we get an
 embedding $M\to N$, which is $S$-close to $\varphi$.
 The restriction $K\to N$ of the embedding is $S$-close to $\varphi$,
 and hence $\varphi$ is approximable by embeddings.
 \end{proof}

 For the proof of Lemma~4.2.I we need the following definition.

 \begin{definition*}[Definition~of~$\varphi^s$]
 For a nondegenerate simplicial map $\varphi:K\to G$ denote by
 $\tilde\Delta(\varphi)=\{\, (x,y)\in \tilde K \,|\, \varphi x=\varphi y \,\}$
 the {\it singular graph}, where $\tilde K$ is the deleted product
 of the graph $K$.
 The vertices of $\tilde\Delta(\varphi)$ are the pairs $(x,y)$
 of vertices of $K$ such that $\varphi x=\varphi y$.
 For an arbitrary simplicial map $\varphi:K\to G$ define
 the simplicial {\it singular} map
 $\varphi^s:\tilde\Delta(\varphi^c)\to G$
 by $\varphi^s(x,y)= \varphi^c x$
 for each vertex $(x,y)\in \tilde\Delta(\varphi^c)$.
 \end{definition*}

 \begin{proof}[Proof of 4.2.I]
 First note that if $\varphi$ does not identify triods
 then $\varphi^s$ does not contain a simple triod.
 Second note that $(\varphi')^s=(\varphi^s)'$
 for any simplicial map $\varphi$ (formally, it follows from
 [5, Proposition~2.11] for the pair of maps $\varphi,\psi$, where $\psi$ is
 the projection $\psi:\tilde K\to K$).
 Third note that $|\tilde\Delta(\varphi^c)|\le k$ and
 $\varphi^s$ is a standard winding for a standard winding $\varphi$.
 Therefore, by Lemma~4.2.T, $(\varphi^{(k)})^s$ is a standard winding,
 and hence $\varphi^{(k)}\Sigma(\varphi^{(k)})$
 is a disjoint union of circles.
 Now note that $\varphi^{(k)}\left|_{\Sigma(\varphi^{(k)})}\right.$ is
 ultra-nondegenerate (because $(\varphi^{(k)})^s$ is not
 ultra-nondegenerate in the opposite case) and
 $\Sigma(\varphi^{(k)})$ has no hanging vertices
 (because $\tilde\Delta(\varphi)$ has a hanging vertex
 in the opposite case).
 So $\varphi^{(k)}\left|_{\Sigma(\varphi^{(k)})}\right.$
 is an ultra-nondegenerate map of a disjoint union of circles into
 a disjoint union of circles $\varphi\Sigma(\varphi)$,
 consequently $\varphi^{(k)}\left|_{\Sigma(\varphi^{(k)})}\right.$
 is a standard winding.
 \end{proof}

 \subsection*{Acknowledgements}

 The author is grateful to Arkady Skopenkov for
 permanent attention to this work.

\end{document}


\title{Об аппроксимируемости вложениями
 циклов на плоскости}
 \author{Михаил Скопенков}
 \address{Department of Differential Geometry, Faculty of Mechanics and Mathematics,
 Moscow State University, 119992, Moscow, Russia,
 and Independent University of Moscow, B. Vlasyevsky, 11, 119002, Moscow, Russia.}
\curraddr{}
\email{skopenkov@rambler.ru}
\thanks{Автор частично поддержан грантом ИНТАС 06-1000014-6277,
грантами Российского Фонда Фундаментальных Исследований 05-01-00993-a,
06-01-72551-НЦНИЛ-а, 07-01-00648-a, Грантом Президента Российской Федерации для государственной поддержки ведущих научных школ Российской Федерации,
проект НШ-4578.2006.1, программой Министерства Образования и Науки ''Развитие научного потенциала высшей школы'',
проект РНП 2.1.1.7988, Фондом поддержки молодых ученых ''Конкурс Мёбиуса''.}
 \subjclass{57Q35 (Primary); 54C25, 57M20 (Secondary)}
 \keywords{Аппроксимируемость вложениями, препятствие Ван Кампена,
 реберный граф, производная графа, производная симплициального отображения, операция $d$,
 трансверсальное самопересечение, стандартная $d$-намотка, симплициальное отображение,
 утолщение графа}
 \begin{abstract} Мы получаем критерий аппроксимируемости вложениями кусочно линейных
 отображений $S^1\to\mathbb{R}^2$, аналогичный доказанному Минцем для кусочно линейных отображений $I\to\mathbb{R}^2$.

 {\bf Теорема.}
 {\it Пусть $\varphi:S^1\to \mathbb{R}^2$ --- кусочно линейное отображение, которое является симплициальным
 для некоторой триангуляции $S^1$ с $k$ вершинами. Отображение $\varphi$
 аппроксимируется вложениями, если и только если для каждого $i=0, \dots, k$ его $i$-я
 производная  $\varphi ^ {(i)} $ (определенная Минцем) не
 содержит трансверсальных самопересечений, и не является стандартной намоткой
 степени $\not\in \{-1,0,1 \} $.}

 Мы выводим из результата Минца
 полноту препятствия Ван Кампена к аппроксимируемости вложениями
 кусочно-линейных отображений $I\to\mathbb{R}^2$.
 \end{abstract}
 \maketitle

 \section{Введение}\label{t-sect1}
 \smallskip

 Кусочно линейное отображение $\varphi:K\to \mathbb{R}^2$ графа $K$ в плоскость
 {\it аппроксимируется вложениями,} если для каждого $\varepsilon> 0 $
 существует отображение $f:K\to\mathbb{R}^2$ без самопересечений,
 $\varepsilon$-близкое к $\varphi$. В большей части статьи мы рассматриваем
 случай, когда $\varphi$ является путем или циклом, то есть, $K\cong I$ или $K\cong S^1$.

 \begin{example}\label{t-ex1-1}
  \cite {11}
 Стандартная $d$-намотка $S^1\to S^1\subset\mathbb{R}^2$
 аппроксимируется вложениями в плоскость, если и только если
 $d\in \{-1,0,1 \} $.

 \end{example}

 Можно доказать также, что симплициальное отображение $S^1\to S^1$
 аппроксимируется вложениями, если и только если его степень $d\in \{-1,0,1 \} $
 (см. Теорему~\ref{t-th1-3}).  {\it Трансверсальным самопересечением} кусочно линейного отображения $\varphi:K\to \mathbb{R}^2$
 называется пара непересекающихся дуг $i, j\subset K$, таких что
 $\varphi i$ и $\varphi j$ пересекаются на плоскости трансверсально.

 \begin{example}\label{t-ex1-2}
 Эйлеров путь или цикл в
 графе на плоскости аппроксимируется вложениями, если и только если он
 не имеет трансверсальных самопересечений
 (следовательно, у любого эйлерова графа на плоскости есть эйлеров цикл,
 аппроксимируемый вложениями).
 \end{example}

 Понятие аппроксимируемости вложениями появилось
 в исследованиях вложимости компактов в $\mathbb{R}^2$
 (см. \cite {11, 14}, актуальные обзоры можно найти в статьях
 \cite{17}, \cite[\S9]{7}, \cite[\S4]{2}, \cite[\S1]{8},
 мы вернемся к этому вопросу еще раз в конце \S\ref{t-sect1}).
 Существует алгоритм проверки того,
 является ли данное симплициальное отображение
 аппроксимируемым вложениями (см. \cite {13}).
 Более удобный для применения критерий аппроксимируемости вложениями
 симплициального пути на плоскости был доказан в статье
 \cite{6} (Теорема~\ref{t-th1-3}.I ниже, обобщающая Пример~\ref{t-ex1-2}). Главный
 результат этой статьи --- аналогичный критерий для аппроксимируемости
 вложениями цикла на плоскости (Теорема~\ref{t-th1-3}.S ниже, также
 обобщающая Пример~\ref{t-ex1-2}). Эти критерии утверждают, что, в некотором смысле,
 трансверсальное самопересечение --- единственное препятствие
 к аппроксимируемости вложениями.
 Ясно, что буквально это не верно \cite {11},
 и нет никакого критерия для рассматриваемой проблемы, аналогичного критерию Куратовского.

\begin {figure}
\includegraphics {skopenkov2.1}
\caption {Производные графов и путей}
\label {t-fig1}
\end {figure}

 Мы формулируем наш критерий (Теорему~\ref{t-th1-3}) в терминах {\it производной}
 пути \cite {Min94, 6} ({\it ''операция $d$''}). Дадим определение этого понятия
 (см. иллюстрацию~\ref{t-fig1}).
 Сначала определим {\it производную} $ G' $ графа $G$
 (это --- синоним для {\it реберного графа}).
 Множество вершин графа $ G' $ находится в 1-1 соответствии
 с множеством ребер графа $G$.
 Для ребра $a\subset G$ обозначим через $a'\in G'$ соответствующую вершину.
 Вершины $a'$ и $b' $ в графе $ G' $
 соединены ребром, если и только если ребра $a$ и $b$ являются смежными
 в графе $G$. Отметим, что производные
 $ G' $ и $H' $ гомеоморфных, но не изоморфных графов $G$
 и $H$ не обязательно гомеоморфны. Обратим внимание, что немного изменим определение графа $G'$ ниже.

 Теперь пусть $\varphi$ --- путь в графе $G$, заданный последовательностью
 своих вершин $x_1, \dots, x_k\in G$, где
 вершины $x_i$ и $x _ {i+1} $ соединены ребром.
 Тогда $ (x_1x_2) ', \dots, (x _ {k-1} x_k) '$ является
 последовательностью вершин графа $ G' $.
 В этой последовательности заменим каждый отрезок вида
 $$ (x_ix _ {i+1}) ', (x _ {i+1} x _ {i+2})', \dots, (x_{j-2}x_{j-1})', (x _ {j-1} x _ {j}) ',$$ где
 $ (x_ix _ {i+1}) '= (x _ {i+1} x _ {i+2})' = \dots = (x _ {j-1} x _ {j}) '$, единственной
 вершиной (то есть заменим несколько идущих подряд одинаковых вершин на одну вершину).
 Полученная последовательность вершин определяет путь в графе $ G' $.
 Этот путь $\varphi '$ называют {\it предпроизводной} пути $\varphi$. Пути $\varphi$ и $\varphi'$ мы будем рассматривать как кусочно линейные отображения отрезка $\varphi\colon I\to G$ и $\varphi'\colon I\to G'$.


 Пусть теперь граф $G$ вложен в плоскость, причем его ребра являются прямолинейными отрезками.
 Производная графа $G$ может и не быть планарным графом: например, производная  5-ода (то есть конуса над 5 точками)
 является непланарным графом Куратовкого.
 Однако, если путь $\varphi$ не имеет трансверсальных самопересечений, то образ отображения $\varphi'\colon I\to G'$ является планарным подграфом
 $ G'_\varphi\subset G '$. Мы строим естественное
 вложение $ G'_\varphi\to\mathbb{R}^2$ в \S\ref{t-sect2}, см. Определение производной утолщения. Заменим граф $ G' $ на образ  $G'_\varphi$, и отображение $\varphi'$ --- на
 ограничение $\varphi '\colon I\to G'_\varphi$, называемое {\it производной} пути $\varphi$. В дальнейшем $ G' $ и $\varphi '$ обозначают именно указанные ограничения.  Определим $k$-ю \emph{производную}
 $\varphi ^ {(k)} $ индуктивно. Разумеется, отображение $\varphi ^ {(k)}\colon I\to\mathbb{R}^2$ определено, только если все предыдущие производные $\varphi, \varphi',\dots, \varphi ^ {(k-1)}\colon I\to\mathbb{R}^2$ не имеют трансверсальных самопересечений.  Для цикла $\varphi$ определение {\it производной} $\varphi '$ аналогично, и это будет снова некоторый цикл в графе на плоскости
 (область определения отображения $\varphi '$ может выродиться в точку или пустое множество).

 Приведем пример, который будет использоваться
 в дальнейшем: $\varphi' = \varphi$ для стандартной
 $d$-намотки $\varphi:S^1 \to S^1$ с $d\ne0$. Ясно, что
 $\varphi '$ --- вложение для любого эйлерового пути или цикла $\varphi$ без трансверсальных самопересечений.
 Таким образом, Пример~\ref{t-ex1-2} --- действительно частный случай следующей теоремы.

 \begin{theorem}\label{t-th1-3}
 I) \cite {6}
 Пусть $\varphi:I\to \mathbb{R}^2$ --- кусочно линейное отображение, являющееся симплициальным для некоторой
 триангуляции отрезка $I$ с $k$ вершинами. Отображение $\varphi$
 аппроксимируется вложениями, если и только если для каждого $i=0, \dots, k$
 его $i$-я производная $\varphi ^ {(i)} $ не содержит
 трансверсальных самопересечений.

 S) Пусть $\varphi:S^1\to \mathbb{R}^2$ --- кусочно линейное отображение, являющееся симплициальным для некоторой
 триангуляции окружности $S^1$ с $k$ вершинами. Отображение $\varphi$
 аппроксимируется вложениями, если и только если для каждого
 $i=0, \dots, k$ его $i$-я производная $\varphi ^ {(i)} $
 не содержит трансверсальных самопересечений,
 и при этом не является стандартной намоткой степени $d\not\in \{-1,0,1 \} $.
 \end{theorem}

 Прокомментируем эту теорему. Отображения $\varphi, \varphi',\varphi'',\dots \colon I\to\mathbb{R}^2$ или $S^1\to\mathbb{R}^2$ корректно определены как отображения в плоскость вплоть до первого номера, на котором у производной обнаружилось трансверсальное самопересечение. Если такой номер нашелся, то проверяемое условие (после слов ``если и только если'') уже не выполнено, и $\varphi$ не аппроксимируется вложениями. Иначе, в теореме~\ref{t-th1-3}.I, проверяемое условие выполнено, и $\varphi$ аппроксимируется вложениями (а в теореме~\ref{t-th1-3}.S нужно еще найти степень намотки $\varphi ^ {(k)} $; ср.~с леммой~\ref{t-lem2-3}).

 Мы доказываем обе теоремы~\ref{t-th1-3}.I и \ref{t-th1-3}.S в \S\ref{t-sect2}.
 Наше доказательство результата~\ref{t-th1-3}.I является более простым, чем приведенное в \cite {6}.

 В \S\ref{t-sect3} мы применяем Теорему~\ref{t-th1-3} для получения следующего критерия.

 \begin{corollary}\label{t-cor1-4}
 Кусочно линейное отображение $\varphi:I\to \mathbb{R}^2$ аппроксимируется
 вложениями, если и только если выполнено любое из следующих эквивалентных
 условий:

 D) (свойство взрезанного произведения)
 Существует отображение $$\{(x, y) \in I\times I: x\ne y \}\to S^1,$$ такое
 что его ограничение на множество
 $\{(x, y) \in I\times I: \varphi x\ne\varphi y \} $
 гомотопно отображению, заданному формулой
 $\tilde\varphi (x, y) = \frac {\varphi x-\varphi y} {\| \varphi x-\varphi y \|} $;

 V) препятствие ван Кампена (определенное в \S\ref{t-sect3}) $v (\varphi) =0$.
 \end{corollary}

 Хотя Критерий \ref{t-cor1-4}.V и труднее сформулировать, но его легче
 применять, чем \ref{t-th1-3}.I и \ref{t-cor1-4}.D.
 В Следствии~\ref{t-cor1-4} отрезок $I$ нельзя заменить на окружность $S^1$: стандартная
 3-намотка является контрпримером \cite {8}.
 Препятствия, подобные \ref{t-cor1-4}.D и \ref{t-cor1-4}.V, существуют и в близкой теории
 аппроксимируемости сингулярными зацеплениями (то есть, отображениями с
 непересекающимися образами связных компонент), но критерии, аналогичные \ref{t-cor1-4}.I и \ref{t-cor1-4}.DV
 для них не верны (Пример~\ref{t-ex3-3} ниже).

\begin{conjecture*}
Кусочно-линейный путь $\varphi:I\to \mathbb{R}^2$ аппроксимируется вложениями, если и только если
 для любой пары дуг $I_1,I_2\subset I$, такой что $I_1\cap I_2=\emptyset$, пара ограничений
 $\varphi:I_1\to \mathbb{R}^2$ и $\varphi:I_2\to \mathbb{R}^2$ аппроксимируется сингулярными зацеплениями.
\end{conjecture*}

Интересно обобщить критерии~\ref{t-th1-3} и~\ref{t-cor1-4} на кусочно линейные отображения $\varphi:K\to G\subset\mathbb{R}^2$,
где $K$ --- произвольный граф (см. частный случай в \cite{Sko03'}).

\begin{conjecture*} Пусть $\varphi:K\to G\subset\mathbb{R}^2$ --- кусочно линейное отображение, являющееся симплициальным относительно некоторой триангуляции графа $K$ с $k$ вершинами. Тогда отображение $\varphi$ аппроксимируется вложениями, если и только если
$v(\varphi)=0$ и для каждого  $i=0, \dots, k$ его $i$-я производная $\varphi ^ {(i)} $ (определенная в \S\ref{t-sect2})
не содержит стандартных намоток степени $d\not\in \{-1,0,1 \} $.
\end{conjecture*}

Если данная гипотеза верна, то {\it кусочно-линейное отображение $\varphi:K\to \mathbb{R}^2$ дерева $K$ аппроксимируется вложениями, если и только если $v(\varphi)=0$} \cite[Problem~4.5]{2}.

\smallskip


 %
 %
 %


 Завершим \S\ref{t-sect1} несколькими замечаниями по поводу истории возникновения
 понятия аппроксимируемости вложениями.
 Дадим  определение разложения 1-мерного компакта
 в {\it обратный предел} и покажем, как понятие
 аппроксимируемости вложениями появляется при исследовании планарности
 этого компакта.
 (Мы не будем использовать это определение в нашей статье.)
 В качестве примера построим {\it 2-адический соленоид Ван Данцига}.
 Возьмем полноторие $T_1\subset\mathbb{R}^3$.
 Пусть $T_2\subset T_1$ --- полноторие, обходящее дважды вдоль оси
 полнотория $T_1$.
 Аналогично, возьмем полноторие $T_3\subset T_2$, обходящее дважды
 вдоль оси полнотория $T_2$.
 Продолжая далее подобным образом, мы получаем бесконечную последовательность
 полноторий $T_1\supset T_2\supset T_3\supset\dots$
 Пересечение всех полноторий $T_i$ является 1-мерным компактом и
 называется {\it 2-адическим соленоидом Ван Данцига}.
 {\it Обратным пределом} бесконечной последовательности
 графов и симплициальных отображений между ними
 $
 \xymatrix@1{
 {K_1} & {K_2}\ar[l]_{\varphi_1} & {K_3} \ar[l]_{\varphi_2} & {\dots} \ar[l]_{\varphi_3}
 }
 $
 мы называем компакт
 $$
 C = \{\, (x_1, x_2, \dots) \in l_2 \,
 : \, x_i\in K_i \text { и } \varphi_ix _ {i+1} =x _ {i} \, \}.
 $$
 Можно видеть из нашего построения, что
 для соленоида Ван Данцига все $K_i\cong S^1$
 и все $\varphi_i$ суть 2-намотки.
 Можно доказать, что любой 1-мерный компакт может быть представлен
 в виде обратного предела.
 Такое представление показывает, что любой 1-мерный компакт может быть вложен
 в $\mathbb{R}^3$. Оно также предоставляет простое достаточное условие планарности данного компакта:
 для каждого положительного целого числа $i$
 должно существовать вложение $f_i:K_i\to \mathbb{R}^2$, такое что отображение
 $f_i\circ \varphi_i$ аппроксимируется вложениями и
 $f _ {i+1} $ является $2^ {-i}$-близким к $f_i\circ\varphi_i$.




\section{Доказательство критерия аппроксимируемости вложениями}\label{t-sect2}

 Теорема~\ref{t-th1-3} следует из Примера~\ref{t-ex1-1}
 и Лемм~\ref{t-lem2-1}, \ref{t-lem2-2} (для $K\cong I, S^1$)
 и~\ref{t-lem2-3}, которые интересны и сами по себе.

 \begin{lemma}\label{t-lem2-1}
 (для $K\cong I$ \cite{6})
 Предположим, что симплициальное отображение $\varphi:K\to G\subset\mathbb{R}^2$
 графа $K\cong S^1$ или $K\cong I$ не имеет
 трансверсальных самопересечений.
 Тогда если $\varphi '$ аппроксимируется вложениями, то и
 $\varphi$ аппроксимируется вложениями.
 \end{lemma}

 \begin{lemma}\label{t-lem2-2}
 A) \cite {6}
 Если симплициальное отображение $\varphi:K\to G\subset\mathbb{R}^2$ аппроксимируется
 вложениями, то и отображение $\varphi '$ аппроксимируется
 вложениями.

 V) Если симплициальное отображение $\varphi:K\to G\subset\mathbb{R}^2$ аппроксимируется
 $\mod 2$-вложениями, то отображение $\varphi '$ аппроксимируется
 $\mod 2$-вложениями.
 \end{lemma}

 Здесь {\it $\mod2$-вложение} --- это отображение общего положения $f:K\to\mathbb{R}^2$,
 такое что для каждой пары $a, b$ непересекающихся ребер графа $K$ множество
 $fa\cap fb$ состоит из четного числа точек.
 Определение  производной для симплициального отображения произвольного графа $K$ приводится ниже.

 \begin{lemma}\label{t-lem2-3}
 Пусть $\varphi:S^1\to G$ --- кусочно линейное отображение,
 которое является симплициальным для некоторой
 триангуляции окружности $S^1$ с $k$ вершинами.
 Тогда либо область определения отображения $\varphi ^ {(k)} $ пуста, либо $\varphi ^ {(k)}$
 является стандартной намоткой степени $d\ne0$.
 \end{lemma}

 Это число $d$ можно рассматривать как обобщение
 {\it степени} для любого симплициального отображения $S^1\to G$.
 Таким образом, интересно получить решение следующей задачи
 (оно может также сделать применение Критериев~\ref{t-th1-3} 
 более удобным):
 {\it
 найти простой алгоритм для вычисления степени
 намотки $\varphi ^ {(\infty)} $ для данного кусочно линейного отображения $\varphi:S^1\to G$}.

 Далее мы используем следующее обобщение определения
 производной (для пути), данного в \S\ref{t-sect1}.

 \begin{definition}[Производная симплициального отображения] \cite {6} (см. иллюстрацию~ \ref{t-fig1}, а также часть
 иллюстрации~\ref{t-fig4} ниже) Пусть дано симплициальное отображение $\varphi:K\to G$. Сначала построим граф $K'_\varphi$, который будет
 областью определения производной $\varphi '$. Под {\it $\varphi$-компонентой}
 графа $K$ мы подразумеваем любую связную компоненту $\alpha$
 множества $\varphi ^ {-1} a$, отображаемую {\it на} $a$, для некоторого ребра $a\subset G$.
 Множество вершин графа $K'_\varphi$ находится в 1-1 соответствии
 с множеством всех $\varphi$-компонент.
 Для $\varphi$-компоненты $\alpha\subset K$ обозначим через
 $\alpha '\in K '_\varphi$ соответствующую вершину.
 Две вершины $\alpha '$ и $\beta' $ соединены ребром в графе $ K'_\varphi$,
 если и только если $\alpha\cap \beta\ne\emptyset$.
 {\it Производная} $\varphi': K '_\varphi\to G' $ --- это
 симплициальное отображение, определенное на вершинах графа $K'_\varphi$
 формулой $\varphi '\alpha' = (\varphi \alpha) '$.
 В дальнейшем заменим $\varphi '$ на сюръективное ограничение
 $\varphi ':K '_\varphi\to \varphi'K '_\varphi$.
 (В оригинальном определении статьи \cite {6} граф $ G' $ обозначается как $D (G) $,
 производная $\varphi '$ как $d [\varphi] $, и граф $K_\varphi' $ как $D (\varphi, K) $.)
 \end{definition}

 \begin{proof} [Доказательство Леммы~\ref{t-lem2-3}]
 Будем говорить, что симплициальное отображение $\varphi:K\to G$ является
 {\it сильно невырожденным}, если для каждого ребра $a\subset K$
 образ $\varphi a$ является ребром $G$ (а не вершиной) и
 для каждой пары $a, b\subset K$
 смежных ребер мы имеем $\varphi a\ne\varphi b$.
 Обозначим через $|K|$ число вершин в графе $K$.
 Ясно, что если $K\cong S^1$,
 то $|K_\varphi' | \le |K|$, причем $|K_\varphi' | = |K|$, только если
 $\varphi$ является сильно невырожденным.
 Поэтому лемму достаточно доказать только в последнем случае
 (потому что случаи, когда $K_\varphi '\cong I$ или граф $K_\varphi '$ является точкой,
 тривиальны).
 В случае сильно невырожденного отображения лемма очевидна, но мы приводим доказательство
 для полноты.

 Докажем, что если сильно невырожденное сюръективное симплициальное отображение $\varphi:K\to G$
 графа $K\cong S^1$ не является стандартной намоткой степени, отличной от нуля,
 то $|G' |> |G|$.
 Заметим, что для сильно невырожденного отображения $\varphi:S^1\to G$ граф $G$
 не содержит висящих вершин.
 Если степень каждой вершины графа $G$ равна двум,
 то $\varphi$ является сильно невырожденным симплициальным отображением $S^1\to S^1$,
 следовательно, $\varphi$ является стандартной намоткой, вопреки нашему предположению.
 Значит, граф $G$ содержит вершину степени по крайней мере 3.
 Тогда, по доказанному выше, число ребер графа $G$ больше
 числа вершин, следовательно, $|G' |> |G|$.
 Поскольку для симплициального на отображения $\varphi:K\to G$ мы имеем
 $1\le |G |\le |K|$, то $|G |, | G' |, \dots, |G ^ {(k)} | \le k$
 (напомним, что по определению отображение $\varphi '$ сюръективно).
 Таким образом, есть только две возможности: любой одна (а значит, и $k$-я тоже)
 из производных  $\varphi, \dots, \varphi ^ {(k)} $ --- стандартная намотка ненулевой степени,
 либо $|G ^ {(k)} | =0$, то есть область определения отображения $\varphi ^ {(k)} $ пуста.
 \end{proof}

\begin {figure}
\includegraphics {skopenkov2.2}
\caption {Производная утолщения графа}
\label {t-fig2}
\end {figure}


 Теперь приведем обещанное в \S\ref{t-sect1} построение вложения
 $ G'_\varphi\to\mathbb{R}^2$. Нам будет удобнее
 рассматривать {\it утолщения} графов, а не вложения
 графов в плоскость. В этом смысле обещанное построение эквивалентно
 построению {\it производной утолщения}
 (см. определение утолщения~$N '$ ниже). Далее мы предполагаем, что фиксировано {\it утолщение
 $N$} графа $G$ на плоскости (то есть, регулярная окрестность
 графа $G\subset\mathbb{R}^2$). Мы также предполагаем, что также фиксировано разложение на ручки
 (обозначаемое через $S$)
 $$
 N = \bigcup\limits _ {x\in\text { множество вершин графа } G} N_x
    \cup\bigcup\limits _ {a\in\text { множество ребер графа } G} N _ {(a)},
 $$
 соответствующее графу $G$,
 где $N_x$ обозначают 2-мерные диски, а
 $N _ {(a)} $ --- присоединенные к ним ленточки.
 Обозначим через
 $N_a$ ограничение $N_x\cup N _ {(a)} \cup N_y$
 утолщения $N$ на ребро $a=xy$.
 Фактически, мы не используем планарность $N$ в последующих рассуждениях.
 Можно считать, что утолщение  $N$ является всего лишь ориентируемым (ориентируемость необходима для
 утверждения Примера~\ref{t-ex1-1}).
 Дадим определение производной $N '$ утолщения $N$.
 Это утолщение $N '$ зависит от симплициального отображения $\varphi:K\to G\subset N$
 и определено корректно, только если $\varphi$ не содержит трансверсальных
 самопересечений. Кроме того,
 в случае произвольного графа $K$ мы должны также предположить,
 что не существует пары дуг $i, j\subset K$
 (не обязательно непересекающихся!), таких что пересечение
 $\varphi i\cap\varphi j$ трансверсально.

 \begin{definition}[Производная утолщения графа]
 (см. иллюстрацию~\ref{t-fig2})
 Пусть $\varphi:K\to G\subset N$ --- симплициальное
 отображение, такое что для любой пары дуг $i, j\subset K$
 пересечение $\varphi i\cap\varphi j$ (возможно пустое)
 не трансверсально.
 Возьмем по диску $N' _ {a'} $
 для каждой вершины $a '\in G'$ и по ленточке $N' _ {(a'b')} $
 для каждого ребра $a'b '\subset G'$.
 Тогда $N '$ вместе с его разложением ручки $S' $
 определяется формулой
 $N' = \bigcup N' _ {a'} \cup\bigcup N' _ {(a'b')} $.
 Здесь мы полагаем $N' _ {(a'b')'} =N _ {(a)} $ для каждого ребра $a\subset G$.
 Для каждой пары $a, b\subset G$ смежных ребер, для которых
 $ (\varphi') ^ {-1} (a'b') \ne\emptyset$,
 мы соединяем два диска $N' _ {a'} $ и $N' _ {b'} $ узкой ленточкой
 $N' _ {(a'b')} $ в $N _ {a\cap b} $.
 Поскольку пересечение дуг $a\cup b$ и
 $c \cup d$ не трансверсально ни для какой пары смежных ребер
 $c, d\subset K$, то мы можем выбрать ленточки $N' _ {(a'b')} $
 так, чтобы они не пересекались для различных ребер $a'b '$.
 \end{definition}

 Это определение можно рассматривать
 как построение вложения $N '\to N$, а также вложения
 $G_\varphi '\to\mathbb{R}^2$.
 Заметим, что разбиение на ручки $S '$ и топологический тип утолщения $N' $ не зависят от
 выбора ленточек $N' _ {(a'b')} $ в нашем определении.
 Альтернативное определение производной $D (N) $ утолщения $N$ из статьи \cite {6}
 не зависит также от выбора отображения $\varphi$.
 Утолщение  $N '$ в нашей статье является подутолщением утолщения $D (N) $ (определение которого приводится в статье \cite {6}),  соответствующим подграфу $G_\varphi '\subset G'$.

 Ясно, что для исследования аппроксимируемости
 вложениями симплициальных отображений $K\to G\subset\mathbb{R}^2$
 достаточно рассмотривать только приближения $f:K\to N$.
 Теперь мы собираемся свести задачу аппроксимируемости
 вложениями данного отображения к задаче существования
 вложения, близкого к нему в некотором смысле ($S$-близкого).

\begin{definition}[$S$-аппроксимация] \cite{6}
 Отображение $f: K \to N$ называется
 {\it $S$-аппроксимацией} отображения $\varphi$,
 или, отображение $f$ {\it $S$-близко} к $\varphi$,
 если выполняются следующие условия:

 \smallskip

 \noindent (i) $fx\subset N _ {\varphi x} $ для каждой вершины или ребра $x$ графа $K$;

 \noindent  (ii) $x\cap f^{-1}N _ {(\varphi x)} $ связно для каждого ребра
       $x$ графа $K$ с невырожденным образом $\varphi x$.
\end{definition}

 Согласно Утверждению~2.9 статьи \cite {6}, отображение
 $\varphi:K\to G$
 аппроксимируется вложениями, если и только если существует вложение
 $f:K\to N$, $S$-близкое к $\varphi$.

 Кусочно линейное отображение $\varphi:K\to N$ называется {\it вырожденным},
 если $\varphi c$ является точкой для некоторого
 ребра $c\subset K$.
 Докажем следующее несложное Утверждение о стягивании ребра~\ref{t-cl2-5},
 которое в некотором смысле позволяет считать, что в Леммах~\ref{t-lem2-1} и~\ref{t-lem2-2}
 отображение $\varphi$ является невырожденным.

\begin {figure}
\begin{tabular}{cc}
\includegraphics {skopenkov2.301} & \includegraphics{skopenkov2.302}\\
a & b
\end{tabular}
\caption {Перестройки вырожденных отображений}
\label {t-fig3}
\end {figure}

\begin{claim}[О стягивании ребра]\label{t-cl2-5}
Пусть $\varphi:K\to G$
 --- симплициальное отображение, такое что $\varphi c$ является точкой для некоторого
 ребра $c\subset K$.
 Пусть $K/c$ --- граф, полученный из графа $K$ стягиванием ребра $c$, и пусть $\varphi/c:K/c\to G$ ---
 соответствующее отображение.
 Тогда

 D) $ K_ {\varphi/c} '=K_\varphi' $, $G_\varphi' =G _ {\varphi/c} '$
 и $ (\varphi/c) '= \varphi' $.

 A) для $K\cong S^1$ или $K\cong I$ отображение $\varphi/c$ аппроксимируется
 вложениями, если и только если $\varphi$ аппроксимируется вложениями.

 K) для произвольного графа $K$, если $\varphi$ аппроксимируется вложениями,
 то $\varphi/c$ аппроксимируется вложениями.

 V) Если $\varphi$ аппроксимируется $\mod2$-вложениями,
 то $\varphi/c$ аппроксимируется $\mod2$-вложениями.
 \end{claim}

 \begin{proof}
 [Доказательство Утверждения~\ref{t-cl2-5}]
 D) очевидно.

 A) Докажем прямую импликацию.
 Пусть $f:K/c\to N$ --- вложение, $S$-близкое к $\varphi/c$.
 Пусть $a\subset K$ --- ребро, смежное с $c$
 (если $c$ --- связная компонента графа $K$, то требуемое утверждение очевидно).
 Добавим новую вершину к ребру $a$ графа $K/c$ (иллюстрация ~ \ref{t-fig3}.a).
 Так как $K\cong S^1$ или $K\cong I$, то
 полученный граф изоморфен $K$
 и вложение $f:K\to N$ --- искомое.
 Обратная импликация --- частный случай утверждения K).

 K) Пусть $f:K\to N$ --- вложение, $S$-близкое к $\varphi$.
 Сделаем перестройку, показанную на иллюстрации ~ \ref{t-fig3}.b.
 Получим вложение $\bar f:K/c\to N$, $S$-близкое к $\varphi/c$.

 V) Пусть $f$ --- $\mod2$-вложение, $S$-близкое к $\varphi$.
 Сделаем перестройку, показанную на иллюстрации~\ref{t-fig3}.b.
 Получим $S$-близкое к $\varphi/c$ отображение $\bar f:K/c\to N$.
 Достаточно доказать, что $ |\bar fa\cap \bar fb | = 0\pmod2$
 для каждой пары непересекающихся ребер $a, b\subset (K/c) $.
 Действительно, $a$ и $b$ являются ребрами также и в графе $K$,
 причем по крайней мере одно из них не
 смежно с $c$ (потому что $a$ и $b$ являются непересекающимися в $K/c$).
 Если ни $a$, ни $b$ не смежно с $c$, то
 $ |\bar fa\cap \bar fb | = | fa\cap fb | = 0\pmod2$.
 Если, например, $b\subset K$ смежно с $c$
 и $a$ не смежно с $c$, то
  $ |\bar f a\cap \bar f b | = | fa\cap fb | + | fa\cap fc | = 0\pmod2$,
 что доказывает утверждение.
 \end{proof}

 Вырожденные отображения появляются в нашем доказательстве Лемм~\ref{t-lem2-1} и~\ref{t-lem2-2},
 даже если исходное отображение $\varphi:K\to G$ является невырожденным.
 Мы собираемся построить
 граф $\bar K '_\varphi$ и пару (вырожденных) симплициальных отображений
 $\xymatrix{ {G} & \bar K '_\varphi \ar[l]_{\bar\varphi} \ar[r]^{\bar\varphi'} & G' }$, которые
 могут быть получены из отображений $\varphi$ и $\varphi '$, соответственно,
 операцией из Утверждения о стягивании ребра~\ref{t-cl2-5}
 (при некоторых дополнительных предположениях относительно $\varphi$,
 детали представлены ниже).
 Вместе с построением вложения $N '\to N$
 (см. определение утолщения~$N '$ выше), это немедленно влечет утверждение Леммы~\ref{t-lem2-1}
 (см. иллюстрации ~ \ref{t-fig4}, \ref{t-fig5}, \ref{t-fig6}).

 \begin{definition}[Полупроизводная симплициального отображения]
 (см. иллюстрацию ~ \ref{t-fig4})
 Предположим, что отображение $\varphi$ является невырожденным, и $K$ не
 имеет вершин степени 0.
 Возьмем несвязное объединение всех $\varphi$-компонент
 графа $K$ (см. Определение производной~$\varphi '$ выше). Соединим ребром любые две вершины,
 принадлежащие различным $\varphi$-компонентам и
 отвечающие одной и той же вершине графа $K$.
 Обозначим полученную {\it полупроизводную} графа $K$ через $\bar K '_\varphi$.
 Таким образом, $\varphi$-компонента $\alpha\subset K$ является также подграфом
 графа $\bar K_\varphi '$, обозначаемым через $\bar\alpha' $.
 В дальнейшем мы отождествляем точки графов $\alpha$ и $\bar\alpha' $.
 Определим симплициальные отображения $\bar\varphi$ и
 $\bar \varphi '$ ({\it полупроизводные} отображения $\varphi$) как очевидные проекции
 $\bar K'_\varphi\to G$ и $\bar K_\varphi '\to G'$, соответственно,
 заданные на вершинах формулами
 $\bar\varphi x =\varphi x$ и
 $\bar \varphi'x = (\varphi \alpha) '$,
 где вершина $x\in\bar K_\varphi '$ принадлежит
 $\varphi$-компоненте $\bar\alpha '$.
 \end{definition}

\begin {figure}
\includegraphics {skopenkov2.4}
\caption {Полупроизводные симплициального отображения}
\label {t-fig4}
\end {figure}


 \begin{proof}[ Доказательство Леммы~\ref{t-lem2-1}]
 Согласно Утверждению о стягивании ребра~\ref{t-cl2-5}.D, отображение $\varphi$ может считать
 невырожденным. Мы также можем считать, что граф
 $K$ не имеет вершин степени 0.
 Легко видеть, что $\varphi$ и $\varphi '$ могут быть получены из
 $\bar\varphi$ и некоторого сужения $\bar\varphi '$,
 соответственно, операцией из Утверждения о стягивании ребра~\ref{t-cl2-5}.
 Если любые две $\varphi$-компоненты имеют не более
 одной общей точки, то $\varphi '$
 может быть получен таким образом непосредственно  из $\bar\varphi$.
 Но для $K\cong S^1$ последнее условие выполнено всегда, кроме случая,
 когда граф $K$ состоит из ровно двух $\varphi$-компонент. Очевидно, отображение $\varphi$
 аппроксимируется вложениями в указанном случае.
 Таким образом, достаточно доказать следующее утверждение:

 (*) если $\bar\varphi '$ аппроксимируется вложениями,
 то $\bar\varphi$ аппроксимируется вложениями.

 Докажем утверждение (*) для произвольного графа $K$.
 Если $\bar \varphi '$ аппроксимируется вложениями,
 то найдется  вложение $\bar K_\varphi '\to N '$, $S '$-близкое к $\bar \varphi' $.
 Определим вложение $f:\bar K_\varphi '\to N$ как композицию
 этого вложения и вложения $N '\to N$,
 построенного в определении утолщения~$N '$ (см. иллюстрацию ~ \ref{t-fig5}, где это построение применяется
 к отображению $\varphi$ с иллюстрации ~ \ref{t-fig4}).
 Ясно, что существует новое разложение на ручки $N=\bigcup \bar N_a \cup \bigcup \bar N_{(ab)}$
 утолщения $N$, обозначаемое $\bar S$, такое что $f$ будет $\bar S$-аппроксимацией отображения
 $\bar\varphi$ (см. иллюстрацию~\ref{t-fig6}, сравни с \cite {6}, Утверждение~4.9)
 Тогда $f:\bar K' _ {\varphi} \to \bar N$ (где $\bar N$ обозначает утолщение
 $N$ с новым разложением ручки $\bar S$) --- вложение,
 $\bar S$-близкое к отображению $\bar\varphi$.
 Лемма доказана.
 \end{proof}


\begin {figure}
\includegraphics {skopenkov2.5}
\caption {Построение $S$-аппроксимации}
\label {t-fig5}
\end {figure}

\begin {figure}
\includegraphics {skopenkov2.6}
\caption {Построение разбиения на ручки}
\label {t-fig6}
\end {figure}

 Та же самая идея используется в доказательстве Лемм~\ref{t-lem2-2}.A,V.
 Рассматривается отображение
 $f:\bar K '_\varphi\to N$ общего положения, $S$-близкое к $\bar\varphi$
 и строится {\it полупроизводная}
 $\bar f ':\bar K '_\varphi\to N '$, $S$-близкая
 к $\bar\varphi '$ (см. иллюстрацию~\ref{t-fig7}).
 Потом проверяется, что если $f$ --- вложение, то
 $\bar f '$ --- также вложение (см. иллюстрацию~\ref{t-fig8}).

 \begin{definition}[Полупроизводная $S$-аппроксимации]
 (см. иллюстрацию ~ \ref{t-fig7},
 где приведенное ниже построение построение применяется к отображению $\varphi$, изображенному на иллюстрации ~ \ref{t-fig4})
 Пусть $K$ --- граф без вершин степени 0.
 Пусть $\varphi:K\to G\subset N$ --- невырожденное
 симплициальное отображение без трансверсальных самопересечений.
 Пусть $f:K\to N$ --- $S$-аппроксимация отображения $\varphi$.
 Тогда {\it полупроизводная} отображения $f$ есть $S '$-аппроксимация
 $\bar f ':\bar K '_\varphi\to N '$ отображения $\varphi'$, и строится следующим образом.
 Для каждого ребра $a\subset G$ выберем гомеоморфизм $h_a:N_a\to N' _ {'} $
 таким образом, что
 для каждого любого ребра $b$, смежного с $a$, мы имеем
 $h_a (N_a\cap N _ {(b)}) \subset N' _ {(a'b')} $.
 Определим $\bar f '$ на каждой $\varphi$-компоненте
 $\bar\alpha '\subset\bar K '_\varphi$ формулой
 $\bar f '\left | _ {\bar\alpha'} \right. =
 h _ {\varphi \alpha} f\left | _ {\alpha} \right.$
 Теперь определим $\bar f '$ на каждом ребре $xy\subset \bar K_\varphi' $,
 соединяющем две различные $\varphi$-компоненты $\bar X '$ и $\bar Y' $.
 Возьмем ребро $a\subset \bar X'$, содержащее вершину $x$.
 Отождествим $\bar X'$ с $X$ (см. Определение полупроизводной симплициального отображения $\bar \varphi' $).
 Тогда $a$ будет отождествлено с некоторым ребром графа $K$, а $x$ --- с некоторой вершиной графа $K$.
 Обозначим через $\bar x$ дугу $a\cap f ^ {-1} N _ {\varphi x} $.
 Определим дугу $\bar y$ аналогично.
 Разрежем ребро $xy$ в три отрезка $xx_1$, $x_1y_1$ и $y_1y$.
 Пусть $\bar f '$ гомеоморфно отображает отрезок $xx_1$ на $h_{\varphi X} f\bar y$,
 отрезок $y_1y$ --- на $h_{\varphi Y} f\bar x$, а отрезок
 $x_1y_1$ --- на прямолинейный отрезок в диске
 $N' _ {(\varphi X \,\varphi Y)}$, соединяющий
 точки $\bar f'x_1$ и $\bar f'y_1$.
 Таким образом, отображение $\bar f ':\bar K_\varphi '\to N '$ построено.
 \end{definition}
 %

\begin {figure}
\begin{tabular}{cc}
\includegraphics {skopenkov2.701} & \includegraphics {skopenkov2.702}
\end{tabular}
\caption {Полупроизводная $S$-аппроксимации}
\label {t-fig7}
\end {figure}


 Заметим, что если $f$ --- вложение, то
 есть более простое альтернативное построение отображения $\bar f '$,
 в некотором смысле обратное к построению из доказательства Леммы~\ref{t-lem2-1}.
 Но это альтернативное построение неприменимо к доказательству Леммы~\ref{t-lem2-2}.V,
 поэтому мы не пользуемся им в данной статье.
 Мы собираемся доказать Лемму~\ref{t-lem2-2}.A,V только в случае,
 когда производная $N '$ определена корректно,
 то есть $K$ не содержит пар дуг $i, j$, для которых пересечение $\varphi i\cap\varphi j$
 трансверсально. Этого достаточно для доказательства Теоремы~\ref{t-th1-3}.
 В общем случае доказательство аналогично, но необходимо всюду вместо $N'$
 пользоваться производной $D (N) $, определенной в статье \cite {6}.

 \begin{proof}[Доказательство Леммы~\ref{t-lem2-2}.A]
 Согласно Утверждению~\ref{t-cl2-5}.K
 можно считать, что $\varphi$ невырождено.
 Возьмем некоторое вложение $f:K\to N$, $S$-близкое к $\varphi$.
 Тогда достаточно показать, что отображение $\bar f '$
 (см. Определение полупроизводной $S$-аппроксимации~$\bar f '$) является вложением.

 Рассмотрим пару различных ребер $xy, zt$ графа $ K'_\varphi$.
 Обозначим множество $\bar f' (xy) \cap \bar f' (zt) $ через $i$.
 Достаточно показать, что $i =\bar f' (xy\cap zt) $.
 Обозначим $a' = \bar\varphi'x$, $b' = \bar\varphi'y$,
 $c' = \bar\varphi'z$ и $d' = \bar\varphi't$.
 Без ограничения общности возможны следующие 3 случая.

 1) $a '$, $b' $, $c '$ и $d' $ попарно различны.
 Так как $\bar f '$ является $S' $-аппроксимацией,то
 $\bar f'xy\subset N' _ {a'b'} $ и
 $\bar f'zt\subset N' _ {c'd'} $, следовательно $i =\emptyset$.

 2) ($a '=c' $ и $b '\ne d '$) или ($a' =b '=c' =d '$).
 Тогда $i\subset N' _ {a'} $, следовательно, $i=h_a (f\bar x\cap f\bar z) $
 (см. определение $h_a$ и $\bar x$ в Определении полупроизводной $S$-аппроксимации ~$\bar f '$,
 дуга $\bar z$ определяется аналогично $\bar x$).
 Если $y\ne t$, то $\bar x$ и $\bar z$ не пересекаются,
 так что $f\bar x\cap f\bar z =\emptyset$ и $i =\emptyset$.
 Если же $y=t$, то $i=h_a (fy) = \bar f' (xy\cap zt) $.

 3) $a '=c' $, $b '=d' $ и $a '\ne b '$.
 В этом случае как $xy$, так и $zt$, соединяют вершины
 из различных $\varphi$-компонент.
 Докажем, что $xy$ и $zt$ не пересекаются.
 Например, пусть $y=t$.
 Тогда все вершины $x$, $y$, $z$ и $t$ графа $\bar K '_\varphi$
 отвечают одной и той же вершине графа $K$. Обозначим ее через $w$.
 Обозначим через $X$ и $Z$ те $\varphi$-компоненты множества
 $\varphi ^ {-1} a=\varphi ^ {-1} c$, для которых $x\in \bar X'$ и $z\in\bar Z'$.
 Таким образом, у $\varphi$-компонент $X$ и $Z$ есть общая точка $w$,
 следовательно $X=Z$.
 Значит, $x, z\in \bar X' = \bar Z' $ отвечают одной и той же вершине $w$,
 следовательно, $x=z$.
 Мы получаем, что $y=t$ и $x=z$, тогда по построению графа $\bar K '_\varphi$
 мы получаем $xy=zt$, что противоречит выбору этих ребер.
 Значит, $xy$ и $zt$ не пересекаются.

 Покажем, что в случае (3) $|i|=0 \pmod2$.
 В дальнейшем будем опускать $\bar f '$ в обозначениях всех образов при отображении $\bar f' $.
 Заметим, что гомеоморфизм $h_a\circ h_b ^ {-1} $
 отображает $y_1y$ и $t_1t$ на $\bar x$ и $\bar z$,
 соответственно (иллюстрация ~ \ref{t-fig8}).
 Из этого следует, что $|i | = | I\cap J|$,
 где $I =\bar x\cup xy_1$ и $J =\bar z\cup zt_1$.
 Из этого также следует, что
 две пары точек $\partial I$ и
 $\partial J$ не зацеплены на окружности $\partial (h_aN _ {a\cap b} \cup N' _ {(a'b')}) $.
 Так как $I, J\subset h_a N _ {a\cap b} \cup N _ {(a'b')} $,
 то $|i|=|I\cap J|=0\pmod2$.
 Таким образом, остается доказать, что $|I\cap J |\le1$, тогда $I\cap J =\emptyset$.
 Последнее утверждение следует из равенства
 $$
 \bar x\cap \bar z=h_a (f\bar x\cap f\bar z) = \emptyset\qquad
 xx_1\cap zz_1=h_a (f\bar y\cap f\bar t) = \emptyset
 \qquad\text {и} \qquad|x_1y_1\cap z_1t_1 |\le1,
 $$
 потому что  $x_1y_1$ и $z_1t_1$ --- прямолинейные отрезки в диске $N _ {(a'b')} $.
 Лемма доказана.
 \end{proof}

\begin {figure}
\includegraphics {skopenkov2.8}
\caption {Подсчет количества точек пересечения}
\label {t-fig8}
\end {figure}


 \begin{proof}[Доказательство Леммы~\ref{t-lem2-2}.V]
 Согласно Утверждению о стягивании ребра~\ref{t-cl2-5}.V
 нам достаточно доказать, что
 если $f:\bar K '_\varphi\to N$ является $\mod2$-вложением, $S$-близким к $\varphi$,
 то его полупроизводная $\bar f '$ также является  $\mod2$-вложением.

 Возьмем пару непересекающихся ребер $xy, zt$ графа $\bar K'_\varphi$
 и рассмотрим те же три случая, что и в доказательстве Леммы~\ref{t-lem2-2}.A.
 Случай 1) тривиален.
 В случае 2) мы имеем $f (xy) \cap f (zt) \subset N_a$, следовательно,
 $|i | = | h_a (f\bar x\cap f\bar z) | = |h_a (f (xy) \cap f (zt)) |
 = |f (xy) \cap f (zt) | =0\pmod2$.
 В доказательстве Леммы~\ref{t-lem2-2}.A мы уже показали, что в случае 3) выполнено равенство
 $|i | = 0\pmod2$. Таким образом, Лемма~\ref{t-lem2-2}.V доказана.
 \end{proof}

 \section{Препятствие Ван Кампена}\label{t-sect3}

 Препятствие Ван Кампена было
 придумано Ван Кампеном при
 исследовании вложимости полиэдров в $\mathbb{R}^{2n} $ \cite {2, FKT94, Kam32,
 7, 8}.
 Дадим определение препятствия ван Кампена
 к аппроксимируемости вложениями симплициальных путей.
 Наше построение более наглядно, чем построение препятствия Ван Кампена к вложимости.
 Пусть $\varphi:I\to \mathbb{R}^2$ --- симплициальный путь (на иллюстрации ~ \ref{t-fig9} приведенная ниже
 конструкция применяется к пути, показанному на иллюстрации ~ \ref{t-fig1}).
 Обозначим через $x_1, \dots, x_k$ вершины графа $I$ в порядке их расположения на дуге $I$,
 и обозначим ребро $x_ix_{i+1} $ через $i$.
 Пусть $I^* =\bigcup\limits_{i<j-1} i\times j$
 --- {\it взрезанный квадрат} графа $I$.
 Раскрасим в красный цвет ребра $x_i\times j$, $j\times x_i$, и клетки $i\times j$
 врезанного квадрата $I ^*$, такие что $\varphi x_i\cap\varphi j =\emptyset$,
 $\varphi i\cap\varphi j =\emptyset$. Обозначим через $I^{*\varphi} $
 красное множество.
 Возьмем отображение общего положения $f:I\to\mathbb{R}^2$, достаточно близкое к $\varphi$.
 В каждую клетку  $i\times j$ ''таблицы'' $I^*$ поставим
 число $v_f (i\times j) = |fi\cap fj |\pmod2$.
 Разрежем $I^*$ вдоль красных ребер.
 Пусть $C_1, C_2, \dots, C_n$ --- все компоненты связности полученной фигуры, для которых
 $\partial C_k\cap\partial I^*\subset I^ {*\varphi} $.
 Обозначим $v_f (C_k) = \sum\limits _ {i\times j\subset C_k}
 v_f (i\times j) $.
 {\it Препятствие Ван Кампена} (с $\mathbb{Z}_2$-коэффициентами) для
 аппроксимируемости вложениями --- это вектор
 $v (\varphi) = (\, v_f (C_1), v_f (C_2), \dots, v_f (C_n) \,) $.

\begin {figure}
\includegraphics {skopenkov2.9}
\caption {Препятствие Ван Кампена}
\label {t-fig9}
\end {figure}


 Несложно проверить, что $v (\varphi) $ не зависит от выбора
 отображения $f$ \cite {8},
 таким образом, $v (\varphi) =0$ является необходимым условием для аппроксимируемости
 вложениями.
 Легко проверить, что $v (\varphi) \ne0$ для кусочно линейного пути $\varphi:I\to\mathbb{R}^2$,
 содержащего трансверсальное самопересечение.
 Таким образом, Следствие~1.4.V следует из~\ref{t-th1-3}, \ref{t-lem2-2}.V и \ref{t-cl3-1}.

 \begin{claim}\label{t-cl3-1}
 Препятствие $v (\varphi) =0$, если и только если
 существует $S$-близкое к отображению
 $\varphi$ $\mod2$-вложение общего положения.
  \end{claim}

 \begin{proof}[Доказательство Утверждения~\ref{t-cl3-1}]
 Обратная импликация очевидна.
 Доказательство прямой импликации следует идеям статьи \cite{Kam32}.
 Мы собираемся использовать
 когомологическую формулировку препятствия Ван Кампена
 (детали можно найти в абзаце перед Утверждением~\ref{t-cl3-2} ниже).
 Пусть $f:K\to N$ --- любая $S$-аппроксимация отображения $\varphi$ общего положения.
 'Движение Райдемайстера', показанное на иллюстрации ~ \ref{t-fig10}.a, добавляет к коциклу $v_f$ кограницу
 $\delta [x\times y] $ элементарной коцепи из группы $B^2 (\tilde K) $.
 Так как $v (\varphi) =0$, то с помощью нескольких таких 'шагов' мы
 можем получить отображение $f:K\to N$, для которого $v_f=0$.
 Тогда отображение $f$ необходимо является $\mod2$-вложением, потому что $v_f=0$ означает, что
 $|fa\cap fb | = 0\pmod2$ для любой пары
 непересекающихся ребер $a, b$ графа $K$.
 \end{proof}

 \begin {figure}
\begin{tabular}{cc}
\includegraphics {skopenkov2.101} & \includegraphics {skopenkov2.102}\\
a & b
\end{tabular}
\caption {''Движение Райдемайстера''}
\label {t-fig10}
\end {figure}


 Теперь мы собираемся доказать, что условия~\ref{t-cor1-4}.V и~\ref{t-cor1-4}.D
 эквивалентны (Утверждение~\ref{t-cl3-2}).
 Утверждение~\ref{t-cl3-2} означает только, что
 \ref{t-cor1-4}.D$\implies$\ref{t-cor1-4}.V,
 но это достаточно для доказательства Следствия~\ref{t-cor1-4}.
 Мы доказываем Утверждение~\ref{t-cl3-2} в более общей формулировке, поэтому нам потребуется
 еще несколько определений.

 Пусть $K$ --- произвольный граф. 
 Пусть $\varphi:K\to G\subset\mathbb{R} ^ {2}$ --- симплициальное отображение.
 Обозначим через $\sigma$ и $\tau$ любые ребра графа $K$. 
 {\it Взрезанным квадратом} графа $K$
 мы называем множество
 $\tilde K =\bigcup \{\,\sigma\times\tau:\sigma\cap\tau =\emptyset \, \} $.
 Пусть $ K^ * =\tilde K/\mathbb{Z}_2$ --- фактор построенного полиэдра относительно антиподального $\mathbb{Z}_2$-действия.
 Пусть $\tilde K ^ {\varphi} \subset \tilde K$ --- подмножество, определяемое формулой
 $\tilde K ^ {\varphi} =
   \{\, \sigma\times\tau:\varphi\sigma\cap\varphi\tau =\emptyset \, \} $.
 Пусть $ K ^ {*\varphi} = \tilde K ^ {\varphi}/\mathbb{Z}_2$.
 Для отображения общего положения $f:K\to\mathbb{R}^ {2} $, близкого к отображению $\varphi$,
 определим коцепь
 $v_f\in C^1 (K ^ *, K ^ {*\varphi}; \mathbb{Z}_2) $ формулой
 $v_f (\sigma\times \tau) =f\sigma\cap f\tau\pmod2$.
 Класс $v (\varphi) = [v_f] \in H^1 (K ^ *, K ^ {*\varphi}; \mathbb{Z}_2) $ этой коцепи
 не зависит от отображения $f$ и называется {\it препятствием Ван Кампена}
 к аппроксимируемости вложениями отображения $\varphi$.
 Мы говорим, что отображение $\varphi:K\to G\subset\mathbb{R}^ {2} $ удовлетворяет
 {\it свойству взрезанного квадрата}, если отображение
 $\tilde\varphi:\tilde K ^ {\varphi} \to S ^ {1} $, заданное формулой
 $\tilde\varphi (x, y) = \frac {\varphi x-\varphi y} {\| \varphi x-\varphi y \|} $,
 продолжается до эквивариантного отображения $\tilde K\to S ^ {1} $.
 Очевидно, данное определение свойства взрезанного квадрата
 эквивалентно~\ref{t-cor1-4}.D в случае $K\cong I$.

 \begin{claim}\label{t-cl3-2}
 Если кусочно линейное отображение $\varphi:K\to \mathbb{R}^{2} $ удовлетворяет свойству взрезанного квадрата,
 то препятствие Ван Кампена  $v (\varphi) =0$.
 \end{claim}

\begin{proof}[Доказательство Утверждения~\ref{t-cl3-2}]

 Возьмем отображение общего положения $f:K\to\mathbb{R}^{2} $, близкое к $\varphi$,
 и определим эквивариантное отображение
 $\tilde f:\tilde K ^ {\varphi} \cup\operatorname{sk} ^ {1} \tilde K\to S ^ {1} $
 формулой $\tilde f (x, y) = \frac {fx-fy} {|fx-fy |}$.
 По общему положению получаем, что $\tilde f$ определено корректно.
 Так как отображение $f$ близко к отображению $\varphi$, то
 $\tilde f\left | _ {\tilde K ^ {\varphi}} \right.$ гомотопно  $\tilde \varphi$.
 Очевидно, если $\tilde\varphi$ продолжается до эквивариантного
 отображения $\tilde K\to S ^ {1} $, то
 $\tilde f\left | _ {\tilde K ^ {\varphi}} \right.$ продолжается до эквивариантного
 отображения $\tilde K\to S ^ {1} $.

 Рассмотрим клетку $\sigma\times\tau\subset\tilde K-\tilde K ^ {\varphi} $,
 где $\sigma, \tau\subset K$ являются $1$-мерными клетками.
 Если отображение $\tilde f$ продолжается на клетку $\sigma\times\tau$, то
 $\deg\tilde f\left |_{\partial (\sigma\times\tau)} \right. =0$.
 Можно показать, что
 $$\deg\tilde f\left |_ {\partial (\sigma\times\tau)} \right.
 =f\sigma\cap f\tau=v_f (\sigma\times \tau) \pmod2. $$
 Поэтому если отображение $\tilde f$ продолжается до эквивариантного отображения
 $\tilde K\to S ^ {1} $, то $v_f=0$.

 Теперь пусть $g:\tilde K ^ {\varphi} \cup\operatorname{sk}^ {1} \tilde K\to S ^ {1} $ ---
 эквивариантное отображение, такое что $gx =\tilde fx$ для каждого
 $x\in\tilde K ^ {\varphi} \cup\operatorname{sk}^ {0} \tilde K$.
 Определим коцепь $v_g\in C ^ {2} (K ^ *, K ^ {*\varphi}; \mathbb{Z}_2) $
 формулой
 $v_g (\sigma) = \deg g\left |_{\partial\sigma} \right.$
 для каждой $2$-мерной клетки $\sigma$.
 Пусть $\sigma\subset\tilde K-\tilde K ^ {\varphi} $ ---
 клетка размерности $1$.
 Возьмем несвязное объединение $\sigma\sqcup\sigma '$ двух копий $\sigma$
 и приклеим $\sigma$ к $\sigma '$ по границе $\partial\sigma =\partial\sigma' $.
 Пусть $d _ {\sigma} $ --- отображение полученной $ 1 $-мерной сферы в $S ^ {1} $,
 заданное формулой $d _ {\sigma} x=fx$ для всех $x\in\sigma$ и
 $d _ {\sigma} x=gx$ для всех $x\in\sigma '$.
 Определим коцепь $v _ {fg} \in C ^ {1} (K ^ *, K ^ {*\varphi}; \mathbb{Z}_2) $
 формулой $v _ {fg} (\sigma) = \deg d _ {\sigma} $.
 Тогда, очевидно, $v_g-v_f =\delta v _ {fg} $.

 Полученная формула означает, что когомологический класс $[v_g] $
 не зависит от выбора эквивариантного отображения
 $g:\tilde K ^ {\varphi} \cup\operatorname{sk}^ {1} \tilde K\to S ^ {1} $ и
 совпадает с препятствием Ван Кампена $v (\varphi) $.
 Это доказывает наше утверждение.
\end{proof}


 \begin{example}\label{t-ex3-3}
 (сравни  c \cite {15, 1})
 {Существует пара кусочно линейных путей $\varphi:I\to \mathbb{R}^2$,
 $\psi:I\to \mathbb{R}^2$
 (см. иллюстрацию~\ref{t-fig11}, где изображена пара путей $f$, $g$, близких к данным),
 не аппроксимируемых сингулярными зацеплениями
 (то есть, отображениями с непересекающимися образами),
 удовлетворяющих следующим условиям:

 V) {\it Препятствие Ван Кампена} $v (\varphi, \psi) =0$.

 D) Отображение
 $\Phi: \{\, (x, y) \in I\times I \, | \,\varphi x\ne\psi y \, \}\to S^1$,
 заданное формулой $\Phi (x, y) = \frac {\varphi x-\psi y}
 {\| \varphi x-\psi y \|} $, гомотопически продолжается до отображения $I\times I\to S^1$.

 I) Пара $\varphi '$, $\psi' $ аппроксимируется сингулярными зацеплениями.}
\end{example}


\begin {figure}
\includegraphics {skopenkov2.11}
\caption {Пара отображений, не аппроксимируемая сингулярными зацеплениями}
\label {t-fig11}
\end {figure}


 \begin{proof}[Доказательство Примера~\ref{t-ex3-3}]
 Пусть $ K, L\cong I$ --- графы с вершинами $k_1, \dots, k_5$
 и $l_1, \dots, l_7$, и пусть $G$ --- граф с вершинами $a_1, \dots, a_6$
 и ребрами $a_1a_2$, $a_1a_3$, $a_1a_4$, $a_1a_5$, $a_2a_3$,
 $a_2a_4$ и $a_2a_6$.
 Требуемые симплициальные отображения $\varphi, \psi$ задаются формулами
 $\varphi k_1=a_1$, $\varphi k_2=a_2$, $\varphi k_3=a_3$, $\varphi k_4=a_1$,
 $\varphi k_5=a_2$ и $\psi l_1=a_5$, $\psi l_2=a_1$, $\psi l_3=a_2$,
 $\psi l_4=a_4$, $\psi l_5=a_1$, $\psi l_6=a_2$, $\psi l_7=a_6$.

 Рассмотрим пару $S$-аппроксимаций $f$ и $g$ отображений $\varphi$ и $\psi$,
 соответственно, показанную на иллюстрации~\ref{t-fig11}.
 Легко видеть, что $|fi\cap gj | = 0\pmod2$ для любой пары ребер
 $i\subset K, j\subset L$.
 Это влечет выполнение обоих свойств~\ref{t-ex3-3}.V и~\ref{t-ex3-3}.D (что доказывается аналогично
 доказательству следствия~\ref{t-cor1-4}, см. также Утверждение~\ref{t-cl3-1}).
 Доказательство свойства~\ref{t-ex3-3}.I --- прямое вычисление.

 Докажем, что пара $\varphi, \psi$ не аппроксимируется
 сингулярными зацеплениями. Предположим противоположное утверждение.
 Пусть $ K _ {13}, K _ {35} \subset K$ и $L _ {14}, L _ {47} \subset L$
 --- дуги между точками $k_1$ и $k_3$, $k_3$ и $k_5$,
 $l_1$ и $l_4$, $l_4$ и $l_7$, соответственно.
 Возьмем малую окрестность графа $\varphi K\cup \psi L$
 на плоскости и выберем ее разложение ручки $S$.
 Обозначим через $A_1$, $A_2$ и $A$ диски разложения на ручки $S$,
 соответствующие вершинам $a_1$, $a_2$ и ребру $a_1a_2$,
 соответственно.
 По аналогу Предложения Минца
 (см. абзац после Определения~$S$-аппроксимации в \S\ref{t-sect2})
 найдутся $S$-аппроксимации $f, g$ отображений
 $\varphi$ и $\psi$, соответственно, с непересекающимися образами.
 Так как $fK _ {13} \cap gL=\emptyset$, то
 пары точек $gL _ {14} \cap\partial (A_1\cup A) $
 и $gL _ {47} \cap\partial (A_1\cup A) $ не зацеплены на окружности
 $\partial (A_1\cup A) $.
 Аналогично, $gL _ {14} \cap\partial A_2$ и $gL _ {47} \cap\partial A_2$
 не зацеплены на окружности $\partial A_2$.
 Значит, $gL _ {14} \cap\partial (A_1\cup A_2\cup A) $ и
 $gL _ {47} \cap\partial (A_1\cup A_2\cup A) $
 не зацеплены на окружности $\partial (A_1\cup A_2\cup A) $.
 Тогда $g$ не может быть $S$-апроксимацией отображения $\psi$.
 Полученное противоречие показывает, что $\varphi$ и $\psi$
 не аппроксимируются сингулярными зацеплениями.
 \end{proof}


\subsection*{Благодарности}
Автор благодарен А. Скопенкову за постоянное внимание к данной работе.